\documentclass[reqno]{amsart}
    \usepackage[utf8]{inputenc}

\title{Knot theory and cluster algebras}
\author{Véronique Bazier-Matte}
\address{Département de mathématiques et de statistiques, Université Laval, Québec (QC), G1V 0A6, Canada}
\email{veronique.bazier-matte.1@ulaval.ca}
\author{Ralf Schiffler}
\thanks{The first author was supported by the NSF grant DMS-1802067. The second author was supported by the NSF grants  DMS-1800860 and DMS-2054561. This work was partially supported by a grant from the Simons Foundation.  The authors would like to thank the Isaac Newton Institute for Mathematical Sciences for support and hospitality during the programme Cluster Algebras and Representation Theory when work on this paper was undertaken. This work was supported by: EPSRC Grant Number EP/R014604/1.
}
\address{Department of Mathematics, University of Connecticut, Storrs, CT 06269-1009, USA}
\email{schiffler@math.uconn.edu}

\usepackage{graphicx}
\usepackage{multirow}% http://ctan.org/pkg/multirow
\usepackage{mathrsfs}
\usepackage{enumerate}
\usepackage{amsmath}
\usepackage{amssymb}
\usepackage{amsfonts}
\usepackage{latexsym}
\usepackage{wrapfig}
\usepackage{float}
\usepackage{subcaption}
\usepackage{verbatim}
\usepackage{hyperref}
\usepackage{url}
\usepackage{bm}
\usepackage{color}
\usepackage[dvipsnames]{xcolor}
\usepackage{pgf,tikz,pgfplots, tikz-cd}
\pgfplotsset{compat=1.13}
\usepackage{mathrsfs}
\usetikzlibrary{arrows, knots}
\usepackage[export]{adjustbox}

\usepackage{xy,color}
\xyoption{poly}
\xyoption{2cell}
\xyoption{all}

\newcommand{\df}{\emph}

\renewcommand{\S}{\mathscr{S}}
\newcommand{\s}{\mathfrak{s}}

\newcommand{\x}{\mathbf{x}}

\newcommand{\y}{\mathbf{y}}
\newcommand{\Z}{\mathbb{Z}}
\newcommand{\rad}{\textup{rad}}
\newcommand{\calg}{\mathcal{G}}
\newcommand{\calf}{\mathcal{F}}
\newcommand{\cala}{\mathcal{A}}
\newcommand{\calx}{\mathcal{X}}
\newcommand{\call}{\mathcal{L}}
\def\xx{\mathbf{x}}
\def\yy{\mathbf{y}}

\def\QQ{\mathbb{Q}}
\def\ZZ{\mathbb{Z}}
\def\PP{\mathbb{P}}
\newcommand{\Axyq}{\mathcal{A}(\mathbf{x},\mathbf{y},Q)}
\newcommand{\xyq}{(\mathbf{x},\mathbf{y},Q)}

\newcommand{\ot}{\leftarrow}
\newcommand{\kb}{\Bbbk}

\newcommand{\za}{\alpha}
\newcommand{\zb}{\beta}
\newcommand{\zd}{\delta}
\newcommand{\zD}{\Delta}
\newcommand{\ze}{\epsilon}
\newcommand{\zg}{\gamma}

\newcommand{\zl}{\lambda}
\newcommand{\zs}{\sigma}

\newcommand{\dimprime}{\textup{dim}^\circ}

\newcommand{\figureeightlabel}{
    \coordinate (1) at (0,0);
    \coordinate (2) at (1,0);
    \coordinate (3) at (2,1);
    \coordinate (4) at (3,1);
    \begin{knot}[
    consider self intersections,
    ignore endpoint intersections=false,
    clip width=5]
    \def\x{45}
    \def\y{60}
    \strand [thick] (1)
    to [out=-\x, in=\x-180, "$4$", swap] (2)
    to [out=\x, in=\x-180, "$5$"] (3)
    to [out=\x, in=180-\x, "$6$"] (4)
    to [out=-\x, in=-\x, "$7$"] (2)
    to [out=180-\x, in=\x, "$8$", swap] (1)
    to [out=\x-180, in=\x-180, "$1$", swap] (3.5,0)
    to [out=\x, in=\x] (4)
    to [out=\x-180, in=-\x, "$2$"] (3)
    to [out=180-\x, in=180-\x, "$3$", swap] (1);
    \flipcrossings{8,9,10,11};
    \end{knot}
}
\newcommand{\figureeight}{
    \coordinate (1) at (0,0);
    \coordinate (2) at (1,0);
    \coordinate (3) at (2,1);
    \coordinate (4) at (3,1);
    \begin{knot}[
    consider self intersections,
    ignore endpoint intersections=false,
    clip width=5]
    \def\x{45}
    \def\y{60}
    \strand [thick] (1)
    to [out=-\x, in=\x-180] (2)
    to [out=\x, in=\x-180] (3)
    to [out=\x, in=180-\x] (4)
    to [out=-\x, in=-\x] (2)
    to [out=180-\x, in=\x] (1)
    to [out=\x-180, in=\x-180] (3.5,0)
    to [out=\x, in=\x] (4)
    to [out=\x-180, in=-\x] (3)
    to [out=180-\x, in=180-\x] (1);
    \flipcrossings{8,9,10,11};
    \end{knot}
}
    
\newtheorem{theorem}{Theorem}[section]
\newtheorem*{theorem*}{Theorem}
\newtheorem{lemma}[theorem]{Lemma}
\newtheorem{proposition}[theorem]{Proposition}
\newtheorem{corollary}[theorem]{Corollary}

\theoremstyle{definition}
\newtheorem{definition}[theorem]{Definition}
\newtheorem{example}[theorem]{Example}
\newtheorem{conjecture}[theorem]{Conjecture}
\newtheorem*{conjecture*}{Conjecture}

\theoremstyle{remark}
\newtheorem{remark}[theorem]{Remark}

\usepackage{tabstackengine}
\setstackEOL{\cr}
\setstackgap{L}{12 pt}
\newcommand{\biginddeux}[2]{
    \begingroup \renewcommand*{\arraystretch}{0.5}
    \begin{matrix} #1 \\ #2 \end{matrix}
    \endgroup}
\newcommand{\bigindtrois}[3]{
    \begingroup \renewcommand*{\arraystretch}{0.5}
    \begin{matrix} #1 \\ #2 \\ #3\end{matrix}
    \endgroup}
\newcommand{\bigindquatre}[4]{
    \begingroup \renewcommand*{\arraystretch}{0.5}
    \begin{matrix} #1 \\ #2 \\ #3 \\ #4\end{matrix}
    \endgroup}

\begin{document}
\maketitle
%\tableofcontents

\begin{abstract}
    We establish a connection between knot theory and cluster algebras via representation theory. To every knot diagram (or link diagram), we associate a cluster algebra  by constructing a quiver with potential. The rank of the cluster algebra is $2n$, where $n$ is the number of crossing points in the knot diagram. We then construct $2n$ indecomposable modules $T(i)$ over the Jacobian algebra of the quiver with potential. For each $T(i)$, we show that the submodule lattice is isomorphic to the corresponding lattice of Kauffman states. We then give a realization of the Alexander polynomial of the knot as a specialization of the $F$-polynomial of $T(i)$, for every $i$. Furthermore, we conjecture that the collection of the $T(i)$ forms a cluster in the cluster algebra whose quiver is isomorphic to the opposite of the initial quiver, and that the resulting cluster automorphism is of order two.
\end{abstract}

\section{Introduction}

We establish a connection from cluster algebras and representation theory to knot theory. Let $K$ be a knot diagram (or link diagram) with $n$ crossings.
A segment of $K$ is a segment of the strand from one crossing point to the next.
We associate to  $K $ a quiver $Q$ with $2n$ vertices, one for each segment of $K$, as well as a potential $W$. 
We then construct $2n$ indecomposable representations $T(i)$, each of which encodes the Alexander polynomial of the link.

To be more precise, each crossing point $p$ of the diagram $K$ gives rise to an oriented cycle $w_p$ of length four in the quiver $Q$ and each region $R$ in $K$ gives rise to an oriented cycle $w_R$ in $Q$ whose length is the number of segments at the region $R$. 
The potential $W$ is the difference of the sum of the crossing point cycles and the sum of the region cycles.

Denote by $B$ the Jacobian algebra of the quiver with potential over an algebraically closed field. Then $B$ is a non-commutative algebra which may be infinite dimensional. 
The representations $T(i)$ are finite-dimensional $B$-modules. We construct them explicitly as representations of the quiver $Q$ by specifying a finite dimensional vector space at every vertex and a linear map for every arrow in $Q$. This construction is a representation theoretic analogue of the construction of the Kauffman states in \cite{K}. The direct sum $T=\oplus T(i)$ is called the link diagram  module of $K$.

Let $\cala$ be the cluster algebra with principal coefficients of the quiver $Q$ as defined in \cite{FZ4}.  $\cala $ is a commutative algebra with a special combinatorial structure. 
It is defined as a subring of a field of rational functions by constructing a set of generators, the cluster variables, via a recursive method called mutation that is determined by the quiver $Q$. Each mutation step connects two sets of $2n$ cluster variables and these sets are called the clusters of $\cala$. 

The cluster variables are Laurent polynomials in two sets of indeterminates $x_i$ and $y_i$, for $i=1,2,\ldots,2n$, with positive integer coefficients \cite{FZ1, LS4}. Moreover, their specialization, obtained by setting all $x_i$ equal to $1$ is a polynomial, called the $F$-polynomial \cite{FZ4}.

It was shown in \cite{DWZ2} that $F$-polynomials can also be computed from modules over our Jacobian algebra $B$. If $M$ is a $B$-module then its $F$-polynomial is

\[F_M=\sum_{\mathbf{e}}\chi(\textup{Gr}_{\mathbf{e}}(M))\,\mathbf{y}^{\underline{\mathbf{e}}},\]
where the sum runs over all dimension vectors $\mathbf{e}=(e_i)_{i=1,2,\dots,2n}$ of submodules of $M$ and 
$\mathbf{y}^{\underline{\mathbf{e}}}=y_1^{e_1}y_2^{e_2}\dots y_{2n}^{e_{2n}}$.
Moreover, $\textup{Gr}_{\mathbf{e}}(M)$ is the quiver Grassmannian of $M$, meaning the variety of all submodules of $M$ whose dimension vector is $\mathbf{e}$, and $\chi $ denotes the Euler characteristic. In general, this Euler characteristic is very hard to compute, because it is known that every projective variety can be realized as a quiver Grassmannian. 

We show that the $F$-polynomials of our $B$-modules $T(i)$ have a much simpler formula, since every submodule is uniquely determined by its dimension vector. Therefore
\[F_{T(i)} = \sum_{L\subset T(i)} \mathbf{y}^{\underline{\dim}\,L},\]
where the sum runs over all submodules $L$ of $T(i)$. 
We write $F_{T(i)}|_t$ for the specialization of the $F$-polynomial at 
\begin{equation}
\label{specializationintro}
y_j=\left\{
\begin{array}
 {ll}
 -t &\textup{if segment $j$ runs from an undercrossing to an overcrossing;}\\  
 -t^{-1} &\textup{if segment $j$ runs from an overcrossing to an undercrossing;}\\
 -1 &\textup{if segment $j$ connects two overcrossings or two undercrossings.} \end{array}
 \right.
\end{equation}

The Alexander polynomial $\zD_K$ of an oriented link diagram $K$ is an important polynomial invariant of the link. 
It is a Laurent polynomial in one variable $t$ with integer coefficients. Introduced by Alexander in \cite{Alex}, it has several equivalent definitions, see for example \cite{Lick}. 
In this paper, we follow Kauffman's approach that realizes the Alexander polynomial as a state sum \cite{K}. 
More recently, the Alexander polynomial has been generalized in the work of Osv\'ath and Szab\'o \cite{OS}, as well as Rasmussen \cite{Ras}, on knot Floer homology.

\bigskip

We are now ready to state our main result. Recall that a link is prime if it cannot be decomposed as the connected sum of two non-trivial links.

\begin{theorem}
 \label{thm 1}
 Let $K$ be a diagram of a prime link. Then, for every segment $i$ of $K$, the Alexander polynomial of $K$ is equal to the specialized $F$-polynomial of the $B$-module $T(i)$.
 That is 
 \[\zD_K=F_{T(i)}|_t.\]
\end{theorem}

We point our that the quiver $Q$, and hence the algebra $B$ and the cluster algebra $\cala$, is not an invariant of the link, because $Q$ depends on the choice of the diagram $K$.  For one, the number of vertices in $Q$ is equal to the number of segments in $K$,  which is not invariant under Reidemeister moves. 
Moreover, the definition of $Q$ does not take into account the difference between an overpass and an underpass in $K$. This difference is only recovered in the specialization (\ref{specializationintro}).

The key step in the proof is the following result which is of interest in its own right. 

\begin{theorem}
 \label{thm 2}
 The lattice of Kauffman states of $K$ relative to a segment $i$ is isomorphic to the lattice of submodules of the $B$-module $T(i)$.
\end{theorem}

An interesting question is how the different $T(i)$ are related to each other. We conjecture the following. Recall that a $B$-module $T$ is called a tilting module if the projective dimension of $T$ is at most one, $\textup{Ext}^1_B(T,T)=0$ and there is an exact sequence $0\to B\to T^0\to T^1$ with $T^0, T^1$ in the additive closure of $T$. Also recall that, given a quiver $Q$, its opposite quiver $Q^{\textup{op}}$ is obtained by reversing the direction of all arrows.
\begin{conjecture}
\begin{itemize}
\item[(a)]  The module $T=\oplus_i T(i)$ is a $B$-tilting module. 
\item[(b)] The Gabriel quiver of the  endomorphism algebra of $T$ is isomorphic to the  quiver $Q^{\textup{op}}$.
\item[(c)] The $2n$ cluster variables in $\cala$ corresponding to $T$ form a cluster.
\item[(d)] There is a permutation $\sigma$ of order two such that the mapping that sends the initial cluster variable $x_i$ to the cluster variable corresponding to $T(\zs(i))$ is a cluster automorphism of order two in the sense of \cite{ASS}.
\end{itemize}
\end{conjecture}
Evidence for the conjecture has been obtained in previous work by  David Whiting and the second author in \cite{SW}. They considered a very special family of links, namely 2-bridge links whose continued fraction has at most two parameters. For a slightly simpler quiver than our quiver $Q$, they constructed some of the modules $T(i)$  and showed that  their direct sum $\oplus T(i)$ can be completed to a tilting module that satisfies the conditions in the conjecture.

As an application, we use  a well-known property of the Alexander polynomial to show the following result that is related to the rank-unimodality conjecture of \cite{MO}. 

\begin{theorem}
 \label{thm 3}
  Let $M$ be a module of Dynkin type $\mathbb{A}_n$ and $\mathcal{L}$ the submodule lattice of $M$. Then
  \[\sum_{L\in \mathcal{L}} (-1)^{h(L)}=\left\{
\begin{array}
 {ll} \pm 1&\textup{if $|\mathcal{L}| $ is odd;}\\
 0&\textup{if $|\mathcal{L}| $ is even,}
\end{array}\right.\]
where  $h(L)= \dim L=\sum_{j\in Q_0} \dim L_j$ is the total dimension of the submodule $L$.
\end{theorem}

\subsubsection*{Relation to other work}
A first connection between cluster algebras and knot theory was given by Kyungyong Lee and the second author in \cite{LS6} in the special case of 2-bridge links. The authors realized another invariant, the Jones polynomial, as a specialization of a cluster variable in a cluster algebra of Dynkin type $\mathbb{A}$. This result was based on an ad hoc construction using the fact that both the 2-bridge links and the cluster variables of type $\mathbb{A}$ are parametrized by continued fractions. We now can see this correspondence as a special case of our general construction, as explained in section~\ref{sect 2bridge}. This provides a more conceptual explanation for the results in \cite{LS6}. However, we do not know how to generalize the Jones polynomial specialization to arbitrary links.

Nagai and Terashima used ancestral triangles constructed from continued fractions to give a formula for the cluster variables of type $\mathbb{A}$ and then defined a specialization that produces the Alexander polynomial of the corresponding 2-bridge link, see \cite{NT}. Our specialization is a generalization of theirs.

 In \cite{CDR}, Cohen, Dasbach and Russel gave a realization of the Alexander polynomial for arbitrary knots as a sum over perfect matchings of the bipartite graph whose vertices are given by the crossing points and the regions of the diagram.  Their graph can be recovered from our quiver  by the methods used for plabic graphs, see for example \cite{FWZ}. In the case of 2-bridge knots, their graph is equivalent to the snake graph associated to the continued fraction in \cite{CS4} and in that case their formula seems to be a special case of the cluster variable expansion formula of \cite{MS} and therefore may be related to ours as well. However, in their approach, the weight of a perfect matching is given by edge weights, which in the cluster algebra setup corresponds to $x$-variables, whereas we use the $y$-variables instead. For arbitrary knots, it is unclear if their formula corresponds to a cluster variable.

All of the articles above  consider a single segment of the link to produce a formula for the invariant. In our approach, we rather aim at a conceptual understanding of the collection of the $2n$ objects given by all of the segments of the link inside the cluster algebra and in the module category of the Jacobian algebra.
\bigskip

The paper is organized as follows. After fixing the notation and recalling certain facts and terminology in section \ref{sect prelim}, we review Kauffman's construction of the state poset and the state polynomial in section~\ref{sect Kauffman states}. In section~\ref{sect B}, we define our quiver with potential and its Jacobian algebra $B$. The link module $T=\oplus_i T(i)$ is constructed in section~\ref{sect knot module}. Section~\ref{sect 3} is devoted to the proof of the lattice isomorphism in Theorem~\ref{thm 2}. Then Theorem~\ref{thm 1} is proved in section~\ref{sect main}. We end the paper with the special case of 2-bridge links and the proof of Theorem~\ref{thm 3} in section~\ref{sect 2bridge}.

\subsection*{Acknowledgments} We thank the anonymous referee, as well as Dylan Rupel, for their useful comments. 

\section{Preliminaries}\label{sect prelim}
We recall basic notions and results from knot theory and cluster algebras.
\subsection{Knots and links}\label{sect knots and links}

 A \emph{knot} is a subset of $\mathbb{R}^3$ that is homeomorphic to a circle.
A \emph{link with $r$ components} is a subset  of $\mathbb{R}^3$  that is homeomorphic to $r$ disjoint circles.   Thus a knot is a link with one component. Links are considered up to ambient isotopy. 
A  link is said to be \emph{prime} if it is not the connected sum of two nontrivial links. 

A \emph{link diagram} $K$ is a projection of the link into the plane, that is injective  except for a finite number of  double points that are called \emph{crossing points}. In addition, the diagram carries the information at each crossing point which of the two strands is on top and which is below. 
A diagram is called \emph{alternating} if traveling along a strand alternates between overcrossings and undercrossings.  A link is called \emph{alternating} if it has an alternating diagram.
A link is said to be  \emph{oriented} if for each component a direction of traveling along the strand is fixed. 

 A \emph{curl} is a monogon in the diagram. We usually assume without loss of generality that our link diagrams are without curls, because one can always remove them (by a Reidemeister I move) without changing the link.

Throughout this paper, we assume that all links are prime and all link diagrams have a finite number of crossing points.

\subsubsection{The Alexander polynomial}\label{sect Alexander}
The Alexander polynomial $\zD$ of an oriented link is a polynomial invariant of the link $\zD\in\Z[t^{\pm \frac{1}{2}}]$ that can be defined in terms of homology, see \cite[Chapter 6]{Lick}. For the original definition of Alexander, see \cite{Alex}.
The Alexander polynomial is defined up to multiplication by a signed power of $t$. 
 
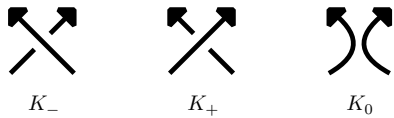
\begin{figure}
\begin{center}
  \scalebox{0.8}{%% Creator: Inkscape 1.0 (4035a4f, 2020-05-01), www.inkscape.org
%% Accompanies image file 'crossing5.pdf' (pdf, eps, ps)
\begingroup%
  \makeatletter%
  \providecommand\color[2][]{%
    \errmessage{(Inkscape) Color is used for the text in Inkscape, but the package 'color.sty' is not loaded}%
    \renewcommand\color[2][]{}%
  }%
  \providecommand\transparent[1]{%
    \errmessage{(Inkscape) Transparency is used (non-zero) for the text in Inkscape, but the package 'transparent.sty' is not loaded}%
    \renewcommand\transparent[1]{}%
  }%
  \providecommand\rotatebox[2]{#2}%
  \newcommand*\fsize{\dimexpr\f@size pt\relax}%
  \newcommand*\lineheight[1]{\fontsize{\fsize}{#1\fsize}\selectfont}%
  \ifx\svgwidth\undefined%
    \setlength{\unitlength}{183.63237762bp}%
    \ifx\svgscale\undefined%
      \relax%
    \else%
      \setlength{\unitlength}{\unitlength * \real{\svgscale}}%
    \fi%
  \else%
    \setlength{\unitlength}{\svgwidth}%
  \fi%
  \global\let\svgwidth\undefined%
  \global\let\svgscale\undefined%
  \makeatother%
  \begin{picture}(1,0.27413715)%
    \lineheight{1}%
    \setlength\tabcolsep{0pt}%
    \put(0,0){\includegraphics[width=\unitlength,page=1]{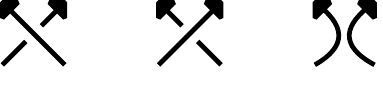}}%
    \put(0.05164099,0.00924141){\color[rgb]{0,0,0}\makebox(0,0)[lt]{\lineheight{0}\smash{\begin{tabular}[t]{l}$K_-$\end{tabular}}}}%
    \put(0.46006573,0.00924141){\color[rgb]{0,0,0}\makebox(0,0)[lt]{\lineheight{0}\smash{\begin{tabular}[t]{l}$K_+$\end{tabular}}}}%
    \put(0.86849044,0.00924141){\color[rgb]{0,0,0}\makebox(0,0)[lt]{\lineheight{0}\smash{\begin{tabular}[t]{l}$K_0$\end{tabular}}}}%
  \end{picture}%
\endgroup%
}
 \caption{Skein relations for the Alexander polynomial.}
 \label{crossing5}
\end{center}
\end{figure}

 In \cite{Conw}, Conway showed that the Alexander polynomial $\zD(K)$ of an oriented link $K$ can be defined recursively 
as follows. The Alexander polynomial of the unknot is $1$, and whenever three oriented links $K_-,K_+$ and $K_0$ are the same except in the neighborhood of a point, where they are as shown in Figure \ref{crossing5}, then 
\[\zD_{K_+}-\zD_{K_-} = (t^{1/2}-t^{-1/2})\zD_{K_0}.\]
This property also provides a normalization of the  Alexander polynomial, but we will not use it here.

The Alexander polynomial $\zD$ has the following properties, see for example \cite[Chapter 6]{Lick}.
\begin{itemize}
\item [(i)] For any link, $\zD(t) \doteq \zD(t^{-1})$, where the symbol $\doteq$ means ``equal up to a signed power of $t$''. 
\item [(ii)]  $\zD(1) = \pm 1$ for any knot, and $\zD(1)=0$ for any link with at least 2 components.
\item [(iii)] For any knot 
\[\zD\doteq a_0+a_1(t^{-1}+t)+a_2(t^{-2}+t^2)+\dots \]
with $a_0$ odd.
\item[(iv)] If a knot has genus $g$ then $2g\ge \textup{breadth}(\zD)$, where the breadth is the difference between the maximal and the minimal degree of the polynomial.
\end{itemize}

Kauffman gave a description of the Alexander polynomial as a state sum. This approach is crucial for us and we review it in section~\ref{sect Kauffman states}.

Let us close this subsection by mentioning a recent breakthrough in a closely related question. In 1982, Freedman showed in \cite{freed} that a knot in the 3-sphere is \emph{topologically} slice if its Alexander polynomial is trivial. 
A famous pair of knots with 11 crossings, the Kinoshita-Terasaka knot and the Conway knot are the smallest non-trivial knots for which the Alexander polynomial is trivial. In particular, both knots are topologically slice. The Kinoshita-Terasaka knot is also known to satisfy the stronger property of being \emph{smoothly} slice.
Recently Lisa Piccirillo solved a longstanding open problem  in \cite{Picc} by proving that the Conway knot is \emph{not} a smoothly slice knot. For an illustration of the quiver of the Conway knot see Example~\ref{ex conway}.

\subsection{Cluster algebras} 
In this section, we recall the definition of a skew-symmetric cluster algebra with principal coefficients following \cite{FZ4, Sln}
  
 Let $\PP$ be the free abelian group on generators $y_1,\ldots, y_n$ written multiplicatively. Let  $\ZZ\PP$ be  the ring of Laurent
polynomials in the variables $y_1,\ldots, y_n$ and let  $\QQ\PP$ denote its field of fractions.
Denote by $\calf=\QQ\PP(x_1,\ldots,x_n)$  the field of rational functions in $n$ variables and coefficients in $\QQ\PP$. We also define an auxiliary addition $\oplus$  by
\begin{equation}
\label{eq:tropical-addition}
\prod_j y_j^{a_j} \oplus \prod_j y_j^{b_j} =
\prod_j y_j^{\min (a_j, b_j)} .
\end{equation}

%The cluster algebra will be defined as a $\ZZ\PP $-subalgebra inside the field $\calf$.
The cluster algebra is determined by the choice of an initial seed $\xyq$, which consists of the following data.
\begin{itemize}
\item $Q$ is a finite connected quiver  without loops  $\xymatrix{\circ\ar@(ur,dr)@<2pt>}$ \hspace{10pt} and 2-cycles $\xymatrix{\circ \ar@<2pt>[r] &\circ \ar@<2pt>[l]}$, and with $n$ vertices;
\item  $\yy=(y_1,\ldots,y_n)$ is the $n$-tuple of generators of $\PP$, called {\em initial coefficient tuple};
\item $\xx=(x_1,\ldots,x_n)$  is the $n$-tuple of variables of $\calf$, called {\em initial cluster}.
\end{itemize}
%The triple $(\xx,\yy, Q)$ is called the \emph{initial seed} of the cluster algebra $\cala=\cala(\xx,\yy, Q)$.

The cluster algebra  $\cala=\cala(\xx,\yy, Q)$
 is the $\ZZ\PP$-subalgebra of $\calf$ generated by so-called \emph{cluster variables}, and these cluster variables are constructed from the initial seed by a recursive method called \emph{mutation}.
A mutation transforms a seed $(\xx,\yy, Q)$ into a new seed $(\xx',\yy', Q')$. Given any seed there are $n$ different mutations $\mu_1,\ldots,\mu_n$, one for each vertex of the quiver, or equivalently, one for each cluster variable in the cluster. %The mutation $\mu_k$ changes the quiver, the coefficient tuple and the cluster of the seed.

The \emph{seed mutation} $\mu_k$ in direction~$k$ transforms
$(\xx, \yy, Q)$ into the seed
$\mu_k(\xx, \yy, Q)=(\xx', \yy', Q')$ defined as follows:

\begin{itemize}
\item $\mathbf{x}'$ is obtained from $\mathbf{x}$ by replacing one cluster variable by a new one, $\mathbf{x}'=\mathbf{x}\setminus\{x_k\}\cup \{x_k'\}$, and $x_k'$ is defined by the following
\emph{exchange relation}

\begin{equation}\label{exch}
{\,x_k} x'_k = {\frac{1}{y_k\oplus 1}\left( y_k\prod_{i\to k} x_i \ 
+ \ \prod_{i\ot k} x_i \right)}
\end{equation}
 where 
 the first product runs over all arrows in $Q$ that end in $k$ and the second product runs over all arrows that start in $k$.
\item $\yy'=(y_1',\dots,y_n')$ is 
a new coefficient $n$-uple, where 
% \begin{equation*}
%\label{eq:y-mutation}
%y'_j =
%\begin{cases}
%y_k^{-1} & \text{if $j = k$};\\[.05in]
%%y_j y_k^{[b_{kj}]_+}
%%(y_k \oplus 1)^{- b_{kj}} & \text{if $j \neq k$}.
%\displaystyle y_j \prod_{k\to j} y_k \frac{\prod_{k\ot j} (y_k\oplus 1)}{\prod_{k\to j} (y_k\oplus 1)}
%& \text{if $j \neq k$}.
%\end{cases}
%\end{equation*}
 \begin{equation*}
\label{eq:y-mutation}
y'_j =
\begin{cases}
y_k^{-1} & \text{if $j = k$};\\[.05in]
%y_j y_k^{[b_{kj}]_+}
%(y_k \oplus 1)^{- b_{kj}} & \text{if $j \neq k$}.
\displaystyle y_j \prod_{k\to j} y_k (y_k\oplus 1)^{-1} \prod_{k\ot j} (y_k\oplus 1)
& \text{if $j \neq k$}.
\end{cases}
\end{equation*}
Note that one of the two products is always trivial, hence equal to 1,  since $Q$ has no oriented 2-cycles. Also note
that $\yy'$  depends only on $\yy$ and $Q$.

\item
The quiver $Q'$ is obtained from $Q$ in three steps: 
\begin{enumerate}
\item for every path $i\to k \to j$ add one arrow $i\to j$,
\item reverse all arrows at $k$,
\item delete 2-cycles.
\end{enumerate} 
\end{itemize}

 Mutations are involutions, that is, $\mu_k\mu_k \xyq=\xyq$. 
Note that  $Q'$ only depends on $Q$,
that $\mathbf{y}'$ depends on $\mathbf{y}$ and $Q$, and that $\mathbf{x}'$ depends on the whole seed $\xyq$.

Let $\calx$ be the set of all cluster variables obtained by finite sequences of mutations from $\xyq$.
Then the cluster algebra $\cala=\Axyq$ is the $\ZZ\PP$-subalgebra of $\calf$ generated by $\calx$.

By definition, the elements of $\cala$ are polynomials in $\calx$ with coefficients in $\ZZ\PP$, so $\cala\subset\ZZ\PP[\calx]$. 
On the other hand, $\cala\subset \calf$, so the elements of $\cala$ are also rational functions in $x_1,\ldots,x_n$ with coefficients in $\QQ\PP$. 
\subsubsection{Laurent phenomenon and positivity}
We have the following important results.
\begin{theorem}[Laurent Phenomenon]\textup{\cite{FZ1}}\label{Laurent}%\textup{[Fomin-Zelevinsky 2002]}
 \, Let $u\in \calx$ be any cluster variable. Then 
\[u=\frac{f(x_1,\ldots,x_n)}{x_1^{d_1}\cdots x_n^{d_n}} \]
where $f\in\ZZ\PP[x_1,\ldots,x_n], d_i\in \ZZ$.\end{theorem}

\begin{theorem}[Positivity]\textup{\cite{LS4}}\label{positivity}%\textup{[Lee-Schiffler 2015]}
\,The coefficients of the Laurent polynomials in Theorem~\ref{Laurent} are positive in the sense that
$f\in\ZZ_{\ge 0}\PP[x_1,\ldots,x_n]$.
\end{theorem}

\subsubsection{$F$-polynomials}
Let $u$ be any cluster variable in the cluster algebra $\cala=\cala\xyq$. By the two theorems above, we can write $u$ as a positive Laurent polynomial in the initial cluster as
$u=\call_u\in\mathbb{Z}_{\ge 0}[x_1,\ldots,x_n,y_1,\ldots,y_n]$. Then the \emph{$F$-polynomial} of $u$ is defined as the evaluation of $\call_u$ at $x_1=\dots=x_n=1$. Thus 
\[F_u=\call_u(1,\ldots,1,y_1,\ldots,y_n).\]

\subsection{Quivers with potential} In this subsection, we recall an alternative approach to $F$-polynomials using quivers with potential.
Let $Q$ be a finite quiver. Following \cite{DWZ1} we let the \emph{vertex span} of $Q$ be the commutative algebra $R$ over $\mathbb{C}$ with basis the constant paths $e_i$, $i\in Q_0$ and multiplication $e_ie_j=\zd_{i,j} e_i$. Furthermore, the \emph{arrow span} of $Q$ is the $R$-bimodule $A$ with $\mathbb{C}$-basis the set of arrows $Q_1$ and $R$-bimodule structure 
$e_iAe_j=\oplus_{\za: j\to i} \mathbb{C}\za$.

The \emph{complete path algebra} of $Q$ then is
$R \langle\langle A\rangle\rangle =\prod_{d=0}^\infty A^{\otimes_R d}$,
with $\mathfrak{m}$-adic topology given by the two-sided ideal $m=\prod_{d=1}^\infty A^{\otimes _R d}$. The elements of $R\langle\langle A\rangle\rangle$ are (possibly infinite) $\mathbb{C}$-linear combinations of paths in $Q$.

A \emph{potential} $W=\sum_{c\in \mathfrak{m}_{\textup{cyc}} }\zl_c c$ on $Q$ is a (possibly infinite) linear combination of cyclic paths in $R\langle\langle A\rangle\rangle$. 
The \emph{cyclic derivative} $\partial_\za$, for $\za\in Q_1$ is %the continuous linear map $\partial_\za\colon\mathfrak{m}_{\textup{cyc}}\to R\langle\langle A\rangle\rangle_{s(\za),t(\za)}$ 
defined on a non-constant cyclic path $\za_1\dots\za_d$ by
\[\partial_\za (\za_1\dots\za_d)=\sum_{p\colon \za_p=\za} \za_{p+1}\dots\za_d\za_1\dots\za_{p-1},\]
and extended linearly to the whole potential.

The \emph{Jacobian algebra} $\textup{Jac}(Q,W)$ of the quiver with potential is defined as the quotient $R\langle\langle A\rangle\rangle / J(W)$,
 where  $J(W)$ is the closure of the two-sided ideal generated by all cyclic derivatives $\partial_\za W$, with $\za\in Q_1$. 
 
 For every finitely generated module $M$ over the Jacobian algebra, 
 Derksen, Weyman and Zelevinsky  introduced  its $F$-polynomial in \cite{DWZ2} as 
 \begin{equation}\label{eq F}
F_M=\sum_{\mathbf{e}} \chi(\textup{Gr}_{\mathbf{e}}(M)) \prod_{i\in Q_0}  y_i^{e_i},
\end{equation}
 where the sum is over all dimension vectors $\mathbf{e}=(e_i)_{i\in Q_0}$ and $\chi(\textup{Gr}_{\mathbf{e}}(M))\in \mathbb{Z}$ is the Euler characteristic of the quiver Grassmannian of all submodules $N \subset M$   of dimension vector $\mathbf{e}$.

 Furthermore, they introduced the notion of mutations of (decorated) representations and showed that if $\mu$ is a mutation  sequence that transforms the zero module into $M$ then the $F$-polynomial of $M$ is equal to the $F$-polynomial of the cluster variable obtained by the same mutation sequence from the initial seed in the cluster algebra $\cala\xyq$.
 
\section{Kauffman states}\label{sect Kauffman states}
 In this section, we recall Kauffman's realization of the Alexander polynomial as a state sum.
%%%%%%%%%%%%%

\subsection{Poset of Kauffman states}\label{sect 1}

Consider an oriented link and fix a  diagram $K$ without curls.
Denote by $n$ the number of crossings. Then, there are $n+2$ regions and $2n$ segments.
We chose a segment and label it $1$ and then, label other segments following the orientation of the string by $2, 3, \dots, 2n$.
In this paper, a pair $(x,R)$ of a crossing point $x$ and a region $R$ such that $x$ is incident to $R$ is called an \df{arrow}.

To define a Kauffman state, we chose a segment $i = 1,2, \dots, n$ and label the adjacent regions $R_i$ and $R_i'$.
A \df{Kauffman state} is a set of arrows $(x,R)$, called \df{markers}, such that:
%A \df{Kauffman state} is a set of \df{markers} where each marker is an arrow $(x,r)$ such that:
\begin{itemize}
    \item each crossing point is used in exactly one marker;
    \item each region except for $R_i$, $R_i'$ is used in exactly one marker.
\end{itemize}
The regions $R_i$, $R_i'$ are used in no marker.

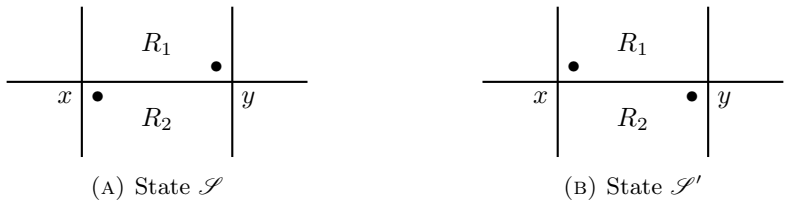
\begin{figure}[ht] 
\centering
\begin{subfigure}{.5\textwidth}
  \centering
  
    \begin{tikzpicture}
    %segments
    \draw[thick] (0,0) -- (4,0);
    \draw[thick] (1,1) -- (1,-1);
    \draw[thick] (3,1) -- (3,-1);
    
    %crossing point
    \coordinate[label=below left:{$x$}] (x) at (1,0);
    \coordinate[label=below right:{$y$}] (y) at (3,0);

    %markers
    \coordinate[label=below right:{$\bullet$}] (m1) at (1,0);
    \coordinate[label=above left:{$\bullet$}] (m2) at (3,0);

%    %markers
%    \coordinate[label=above right:{$\bullet$}] (m1) at (1,0);
%    \coordinate[label=below left:{$\bullet$}] (m2) at (3,0);
    
    %regions
    \node at (2,0.5) {$R_1$};
    \node at (2,-0.5) {$R_2$};
    
    \end{tikzpicture}
  
  %\captionof{figure}{State $\S$}
  \label{fig:state1}
\end{subfigure}%
\begin{subfigure}{.5\textwidth}
  \centering
  
    \begin{tikzpicture}
    %segments
    \draw[thick] (0,0) -- (4,0);
    \draw[thick] (1,1) -- (1,-1);
    \draw[thick] (3,1) -- (3,-1);
    
    %crossing point
    \coordinate[label=below left:{$x$}] (x) at (1,0);
    \coordinate[label=below right:{$y$}] (y) at (3,0);
    
    %markers
    \coordinate[label=above right:{$\bullet$}] (m1) at (1,0);
    \coordinate[label=below left:{$\bullet$}] (m2) at (3,0);
    
    %regions
    \node at (2,0.5) {$R_1$};
    \node at (2,-0.5) {$R_2$};
    \end{tikzpicture}
  
 % \captionof{figure}{State $\S'$}
  \label{fig:state2}
\end{subfigure}
\caption{Kauffman counterclockwise transposition from state $\S$ (left) to state $\S'$ (right).}\label{fig::KauffmanTransposition}
\end{figure}

A state $\S'$ is obtained from a state $\S$ by a \df{counterclockwise transposition} at a segment $j$ if $\S'$ is obtained from $\S$ by switching two markers at the segment $j$ as in Figure \ref{fig::KauffmanTransposition}.

%\[ \]
More precisely, let $x$, $y$ be the endpoints of the segment $j$ and let $R_1$, $R_2$ be the adjacent regions at $j$ such that, going clockwise around $x$, we go from $R_1$ to $R_2$ crossing $j$. Then, $\S$ contains the markers $(x, R_2)$, $(y, R_1)$, $\S'$ contains the markers $(x, R_1)$, $(y, R_2)$ and the other markers in $\S$ and $\S'$ are the same.

We define a partial order on the set of all Kauffman states by $\S_1 < \S_2$ if there is a sequence of counterclockwise transpositions that transforms $\S_1$ into $\S_2$. Kauffman proved that the resulting poset is a lattice whose maximal element is a state that admits no counterclockwise transposition and is therefore called the clocked state in \cite{K}. Similarly, the minimal element is called the counterclocked state in \cite{K}. We will refer to these states as the \emph{maximal} and the \emph{minimal} state.

\begin{example} \label{ex::kstateslattice}
Let's use the following labeling for the segments of the figure-eight knot.
\[\begin{tikzpicture}[scale = 1.2]
\figureeightlabel
\end{tikzpicture}\]

Figure \ref{fig::kstateslattice} shows the lattice of Kauffman states for the figure-eight knot with regards to segment 1.
\begin{figure} 
    \centering
    \begin{tikzpicture}[scale = 2.5]
    \def\sca{0.75}
    \node (n0) at (0,0) {
        \begin{tikzpicture}[scale = \sca]
        \figureeight
        \coordinate[label=above:{$\bullet$}] (m1) at (1);
        \coordinate[label=left:{$\bullet$}] (m2) at (2);
        \coordinate[label=right:{$\bullet$}] (m3) at (3);
        \coordinate[label=below:{$\bullet$}] (m4) at (4);
        \end{tikzpicture}
    };
    \node (n1) at (1,-1) {
        \begin{tikzpicture}[scale = \sca]
        \figureeight
        \coordinate[label=right:{$\bullet$}] (m1) at (1);
        \coordinate[label=above:{$\bullet$}] (m2) at (2);
        \coordinate[label=right:{$\bullet$}] (m3) at (3);
        \coordinate[label=below:{$\bullet$}] (m4) at (4);
        \end{tikzpicture}
    };
    \node (n2) at (-1,-1) {
        \begin{tikzpicture}[scale = \sca]
        \figureeight
        \coordinate[label=above:{$\bullet$}] (m1) at (1);
        \coordinate[label=left:{$\bullet$}] (m2) at (2);
        \coordinate[label=below:{$\bullet$}] (m3) at (3);
        \coordinate[label=left:{$\bullet$}] (m4) at (4);
        \end{tikzpicture}
    };
    \node (n3) at (0,-2) {
        \begin{tikzpicture}[scale = \sca]
        \figureeight
        \coordinate[label=right:{$\bullet$}] (m1) at (1);
        \coordinate[label=above:{$\bullet$}] (m2) at (2);
        \coordinate[label=below:{$\bullet$}] (m3) at (3);
        \coordinate[label=left:{$\bullet$}] (m4) at (4);
        \end{tikzpicture}
    };
    \node (n4) at (0,-3) {
        \begin{tikzpicture}[scale = \sca]
        \figureeight
        \coordinate[label=right:{$\bullet$}] (m1) at (1);
        \coordinate[label=right:{$\bullet$}] (m2) at (2);
        \coordinate[label=left:{$\bullet$}] (m3) at (3);
        \coordinate[label=left:{$\bullet$}] (m4) at (4);
        \end{tikzpicture}
    };
    \draw (n0) -- (n1) -- (n3) -- (n4);
    \draw (n0) -- (n2) -- (n3);
    \end{tikzpicture}
    \caption{Lattice of Kauffman states of the figure-eight knot.}
    \label{fig::kstateslattice}
\end{figure}
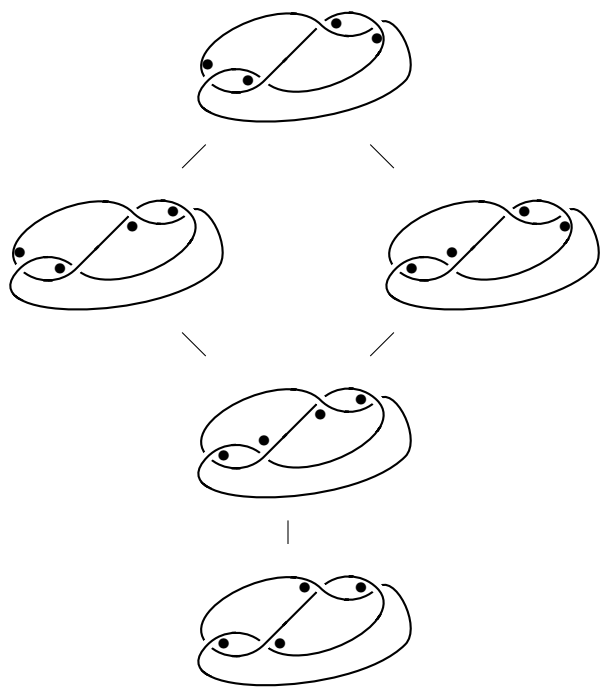

\end{example}

%%%%%%%%%%%%%

\subsection{The State polynomial}

Following Kauffman, we define the \df{weight} $w(x,R)$ of an arrow $(x,R)$ as shown in the following two cases.
\[
\begin{tikzpicture}
%Knot
\begin{knot}[
clip width=5,
flip crossing=1,
]
\strand[thick, -stealth] (-1,-1) -- (1,1);
\strand[thick, -stealth] (1,-1) -- (-1,1);
\end{knot}

%crossing point
\coordinate[label= left:{$x$}] (x) at (0,0);

%regions
\node at (0,1) {$R_1$};
\node at (-1,0) {$R_2$};
\node at (0,-1) {$R_3$};
\node at (1,0) {$R_4$};
\end{tikzpicture}
\]
In this case, $w(x,R_1) = B$, $w(x,R_2) = 1$, $w(x,R_3) = W$ and $w(x,R_4) = 1$.

\[
\begin{tikzpicture}
%Knot
\begin{knot}[
clip width=5,
flip crossing=1,
]
\strand[thick, -stealth] (1,-1) -- (-1,1);
\strand[thick, -stealth] (-1,-1) -- (1,1);
\end{knot}

%crossing point
\coordinate[label= left:{$x$}] (x) at (0,0);

%regions
\node at (0,1) {$R_1$};
\node at (-1,0) {$R_2$};
\node at (0,-1) {$R_3$};
\node at (1,0) {$R_4$};
\end{tikzpicture}
\]
In this case, $w(x,R_1) = W$, $w(x,R_2) = 1$, $w(x,R_3) = B$ and $w(x,R_4) = 1$.

The \df{weight} $w(\S)$ of a state $\S$ is defined as \[ w(\S) = \prod_{(x,R) \in \S} w(x,R).\]
The \df{state polynomial} is the sum of the weights of all states $\S$ \[ \sum_{\S} \sigma(\S) w(\S), \] where $\sigma(\S) = (-1)^b$ with $b$ is the exponent of $B$ in $w(\S)$.

\begin{theorem}[\cite{K}]
The Alexander-Conway polynomial of a diagram $L$ is equal to the specialization of the state polynomial at $W = t^{\frac{1}{2}}$, $B = t^{-\frac{1}{2}}$.
\end{theorem}

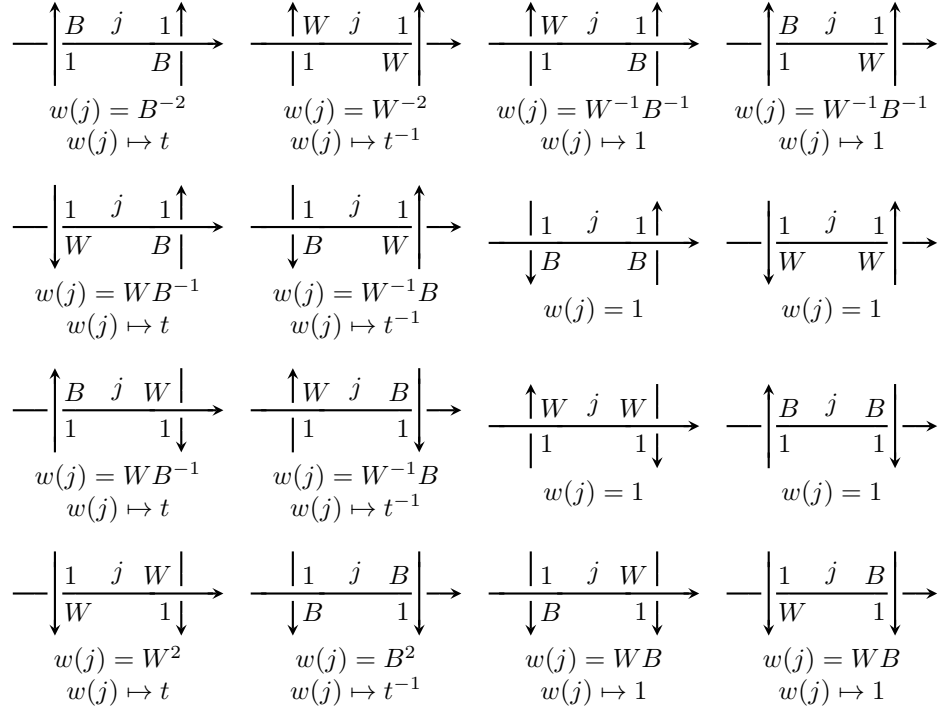
\begin{figure}
\centering
\begin{subfigure}{.25\textwidth}
  \centering
  
    \begin{tikzpicture} [scale=0.7]
    \begin{knot} [clip width=7, flip crossing = 1]
    \strand[thick, -stealth] (0,0) to ["$j$"] (4,0);
    \strand[thick, -stealth] (0.8,-0.8) -- (0.8,0.8);
    \strand[thick, -stealth] (3.2,-0.8) -- (3.2,0.8);
    \end{knot}
    \coordinate[label = above right:{$B$}] (1) at (0.8,0);
    \coordinate[label = below right:{$1$}] (1) at (0.8,0);
    \coordinate[label = below left:{$B$}] (1) at (3.2,0);
    \coordinate[label = above left:{$1$}] (1) at (3.2,0);
    \end{tikzpicture} 
  
  $w(j) = B^{-2}$
  
  $w(j)\mapsto t$
\end{subfigure}%
 \begin{subfigure}{.25\textwidth}
  \centering
  
    \begin{tikzpicture} [scale=0.7]
    \begin{knot} [clip width=7, flip crossing = 2]
    \strand[thick, -stealth] (0,0) to ["$j$"] (4,0);
    \strand[thick, -stealth] (0.8,-0.8) -- (0.8,0.8);
    \strand[thick, -stealth] (3.2,-0.8) -- (3.2,0.8);
    \end{knot}
    \coordinate[label = above right:{$W$}] (1) at (0.8,0);
    \coordinate[label = below right:{$1$}] (1) at (0.8,0);
    \coordinate[label = below left:{$W$}] (1) at (3.2,0);
    \coordinate[label = above left:{$1$}] (1) at (3.2,0);
    \end{tikzpicture} 
  
  $w(j) = W^{-2}$
  
  $w(j)\mapsto t^{-1}$
\end{subfigure}%
\begin{subfigure}{.25\textwidth}
  \centering
  
    \begin{tikzpicture} [scale=0.7]
    \begin{knot} [clip width=7]
    \strand[thick, -stealth] (0,0) to ["$j$"] (4,0);
    \strand[thick, -stealth] (0.8,-0.8) -- (0.8,0.8);
    \strand[thick, -stealth] (3.2,-0.8) -- (3.2,0.8);
    \end{knot}
    \coordinate[label = above right:{$W$}] (1) at (0.8,0);
    \coordinate[label = below right:{$1$}] (1) at (0.8,0);
    \coordinate[label = below left:{$B$}] (1) at (3.2,0);
    \coordinate[label = above left:{$1$}] (1) at (3.2,0);
    \end{tikzpicture} 
  
  $w(j) = W^{-1}B^{-1}$

  $w(j)\mapsto 1$
\end{subfigure}%
\begin{subfigure}{.25\textwidth}
  \centering
  
    \begin{tikzpicture} [scale=0.7]
    \begin{knot} [clip width=7]
    \flipcrossings{1,2}
    \strand[thick, -stealth] (0,0) to ["$j$"] (4,0);
    \strand[thick, -stealth] (0.8,-0.8) -- (0.8,0.8);
    \strand[thick, -stealth] (3.2,-0.8) -- (3.2,0.8);
    \end{knot}
    \coordinate[label = above right:{$B$}] (1) at (0.8,0);
    \coordinate[label = below right:{$1$}] (1) at (0.8,0);
    \coordinate[label = below left:{$W$}] (1) at (3.2,0);
    \coordinate[label = above left:{$1$}] (1) at (3.2,0);
    \end{tikzpicture} 
  
  $w(j) = W^{-1}B^{-1}$

  $w(j)\mapsto 1$
\end{subfigure}%
\vspace{10 pt}
\begin{subfigure}{.25\textwidth}
  \centering
  
    \begin{tikzpicture} [scale=0.7]
    \begin{knot} [clip width=7, flip crossing = 1]
    \strand[thick, -stealth] (0,0) to ["$j$"] (4,0);
    \strand[thick, stealth-] (0.8,-0.8) -- (0.8,0.8);
    \strand[thick, -stealth] (3.2,-0.8) -- (3.2,0.8);
    \end{knot}
    \coordinate[label = above right:{$1$}] (1) at (0.8,0);
    \coordinate[label = below right:{$W$}] (1) at (0.8,0);
    \coordinate[label = below left:{$B$}] (1) at (3.2,0);
    \coordinate[label = above left:{$1$}] (1) at (3.2,0);
    \end{tikzpicture} 
  
  $w(j) = WB^{-1}$
  
  $w(j)\mapsto t$
\end{subfigure}%
 \begin{subfigure}{.25\textwidth}
  \centering
  
    \begin{tikzpicture} [scale=0.7]
    \begin{knot} [clip width=7, flip crossing = 2]
    \strand[thick, -stealth] (0,0) to ["$j$"] (4,0);
    \strand[thick, stealth-] (0.8,-0.8) -- (0.8,0.8);
    \strand[thick, -stealth] (3.2,-0.8) -- (3.2,0.8);
    \end{knot}
    \coordinate[label = above right:{$1$}] (1) at (0.8,0);
    \coordinate[label = below right:{$B$}] (1) at (0.8,0);
    \coordinate[label = below left:{$W$}] (1) at (3.2,0);
    \coordinate[label = above left:{$1$}] (1) at (3.2,0);
    \end{tikzpicture} 
  
  $w(j) = W^{-1}B$
 
  $w(j)\mapsto t^{-1}$
\end{subfigure}%
\begin{subfigure}{.25\textwidth}
  \centering
  
    \begin{tikzpicture} [scale=0.7]
    \begin{knot} [clip width=7]
    \strand[thick, -stealth] (0,0) to ["$j$"] (4,0);
    \strand[thick, stealth-] (0.8,-0.8) -- (0.8,0.8);
    \strand[thick, -stealth] (3.2,-0.8) -- (3.2,0.8);
    \end{knot}
    \coordinate[label = above right:{$1$}] (1) at (0.8,0);
    \coordinate[label = below right:{$B$}] (1) at (0.8,0);
    \coordinate[label = below left:{$B$}] (1) at (3.2,0);
    \coordinate[label = above left:{$1$}] (1) at (3.2,0);
    \end{tikzpicture} 
  
  $w(j) = 1$

\end{subfigure}%
\begin{subfigure}{.25\textwidth}
  \centering
  
    \begin{tikzpicture} [scale=0.7]
    \begin{knot} [clip width=7]
    \flipcrossings{1,2}
    \strand[thick, -stealth] (0,0) to ["$j$"] (4,0);
    \strand[thick, stealth-] (0.8,-0.8) -- (0.8,0.8);
    \strand[thick, -stealth] (3.2,-0.8) -- (3.2,0.8);
    \end{knot}
    \coordinate[label = above right:{$1$}] (1) at (0.8,0);
    \coordinate[label = below right:{$W$}] (1) at (0.8,0);
    \coordinate[label = below left:{$W$}] (1) at (3.2,0);
    \coordinate[label = above left:{$1$}] (1) at (3.2,0);
    \end{tikzpicture} 
  
  $w(j) = 1$

\end{subfigure}%

\vspace{10 pt}

\begin{subfigure}{.25\textwidth}
  \centering
  
    \begin{tikzpicture} [scale=0.7]
    \begin{knot} [clip width=7, flip crossing = 1]
    \strand[thick, -stealth] (0,0) to ["$j$"] (4,0);
    \strand[thick, -stealth] (0.8,-0.8) -- (0.8,0.8);
    \strand[thick, stealth-] (3.2,-0.8) -- (3.2,0.8);
    \end{knot}
    \coordinate[label = above right:{$B$}] (1) at (0.8,0);
    \coordinate[label = below right:{$1$}] (1) at (0.8,0);
    \coordinate[label = below left:{$1$}] (1) at (3.2,0);
    \coordinate[label = above left:{$W$}] (1) at (3.2,0);
    \end{tikzpicture} 
  
  $w(j) = WB^{-1}$
   
  $w(j)\mapsto t$
\end{subfigure}%
 \begin{subfigure}{.25\textwidth}
  \centering
  
    \begin{tikzpicture} [scale=0.7]
    \begin{knot} [clip width=7, flip crossing = 2]
    \strand[thick, -stealth] (0,0) to ["$j$"] (4,0);
    \strand[thick, -stealth] (0.8,-0.8) -- (0.8,0.8);
    \strand[thick, stealth-] (3.2,-0.8) -- (3.2,0.8);
    \end{knot}
    \coordinate[label = above right:{$W$}] (1) at (0.8,0);
    \coordinate[label = below right:{$1$}] (1) at (0.8,0);
    \coordinate[label = below left:{$1$}] (1) at (3.2,0);
    \coordinate[label = above left:{$B$}] (1) at (3.2,0);
    \end{tikzpicture} 
  
  $w(j) = W^{-1}B$
   
  $w(j)\mapsto t^{-1}$
\end{subfigure}%
\begin{subfigure}{.25\textwidth}
  \centering
  
    \begin{tikzpicture} [scale=0.7]
    \begin{knot} [clip width=7]
    \strand[thick, -stealth] (0,0) to ["$j$"] (4,0);
    \strand[thick, -stealth] (0.8,-0.8) -- (0.8,0.8);
    \strand[thick, stealth-] (3.2,-0.8) -- (3.2,0.8);
    \end{knot}
    \coordinate[label = above right:{$W$}] (1) at (0.8,0);
    \coordinate[label = below right:{$1$}] (1) at (0.8,0);
    \coordinate[label = below left:{$1$}] (1) at (3.2,0);
    \coordinate[label = above left:{$W$}] (1) at (3.2,0);
    \end{tikzpicture} 
  
  $w(j) = 1$
\end{subfigure}%
\begin{subfigure}{.25\textwidth}
  \centering
  
    \begin{tikzpicture} [scale=0.7]
    \begin{knot} [clip width=7]
    \flipcrossings{1,2}
    \strand[thick, -stealth] (0,0) to ["$j$"] (4,0);
    \strand[thick, -stealth] (0.8,-0.8) -- (0.8,0.8);
    \strand[thick, stealth-] (3.2,-0.8) -- (3.2,0.8);
    \end{knot}
    \coordinate[label = above right:{$B$}] (1) at (0.8,0);
    \coordinate[label = below right:{$1$}] (1) at (0.8,0);
    \coordinate[label = below left:{$1$}] (1) at (3.2,0);
    \coordinate[label = above left:{$B$}] (1) at (3.2,0);
    \end{tikzpicture} 
  
  $w(j) = 1$
\end{subfigure}%

\vspace{10 pt}

\begin{subfigure}{.25\textwidth}
  \centering
  
    \begin{tikzpicture} [scale=0.7]
    \begin{knot} [clip width=7, flip crossing = 1]
    \strand[thick, -stealth] (0,0) to ["$j$"] (4,0);
    \strand[thick, stealth-] (0.8,-0.8) -- (0.8,0.8);
    \strand[thick, stealth-] (3.2,-0.8) -- (3.2,0.8);
    \end{knot}
    \coordinate[label = above right:{$1$}] (1) at (0.8,0);
    \coordinate[label = below right:{$W$}] (1) at (0.8,0);
    \coordinate[label = below left:{$1$}] (1) at (3.2,0);
    \coordinate[label = above left:{$W$}] (1) at (3.2,0);
    \end{tikzpicture} 
  
  $w(j) = W^2$

  $w(j)\mapsto t$
\end{subfigure}%
 \begin{subfigure}{.25\textwidth}
  \centering
  
    \begin{tikzpicture} [scale=0.7]
    \begin{knot} [clip width=7, flip crossing = 2]
    \strand[thick, -stealth] (0,0) to ["$j$"] (4,0);
    \strand[thick, stealth-] (0.8,-0.8) -- (0.8,0.8);
    \strand[thick, stealth-] (3.2,-0.8) -- (3.2,0.8);
    \end{knot}
    \coordinate[label = above right:{$1$}] (1) at (0.8,0);
    \coordinate[label = below right:{$B$}] (1) at (0.8,0);
    \coordinate[label = below left:{$1$}] (1) at (3.2,0);
    \coordinate[label = above left:{$B$}] (1) at (3.2,0);
    \end{tikzpicture} 
  
  $w(j) = B^2$

  $w(j)\mapsto t^{-1}$
\end{subfigure}%
\begin{subfigure}{.25\textwidth}
  \centering
  
    \begin{tikzpicture} [scale=0.7]
    \begin{knot} [clip width=7]
    \strand[thick, -stealth] (0,0) to ["$j$"] (4,0);
    \strand[thick, stealth-] (0.8,-0.8) -- (0.8,0.8);
    \strand[thick, stealth-] (3.2,-0.8) -- (3.2,0.8);
    \end{knot}
    \coordinate[label = above right:{$1$}] (1) at (0.8,0);
    \coordinate[label = below right:{$B$}] (1) at (0.8,0);
    \coordinate[label = below left:{$1$}] (1) at (3.2,0);
    \coordinate[label = above left:{$W$}] (1) at (3.2,0);
    \end{tikzpicture} 
  
  $w(j) = WB$

  $w(j)\mapsto 1$
\end{subfigure}%
\begin{subfigure}{.25\textwidth}
  \centering
  
    \begin{tikzpicture} [scale=0.7]
    \begin{knot} [clip width=7]
    \flipcrossings{1,2}
    \strand[thick, -stealth] (0,0) to ["$j$"] (4,0);
    \strand[thick, stealth-] (0.8,-0.8) -- (0.8,0.8);
    \strand[thick, stealth-] (3.2,-0.8) -- (3.2,0.8);
    \end{knot}
    \coordinate[label = above right:{$1$}] (1) at (0.8,0);
    \coordinate[label = below right:{$W$}] (1) at (0.8,0);
    \coordinate[label = below left:{$1$}] (1) at (3.2,0);
    \coordinate[label = above left:{$B$}] (1) at (3.2,0);
    \end{tikzpicture} 
  
  $w(j) = WB$
   
  $w(j)\mapsto 1$
\end{subfigure}
\caption{Possible values and specializations at $W = t^{\frac{1}{2}}$ and $B = t^{-\frac{1}{2}}$ for the weight of a segment $j$.}
\label{fig::weightsegment}
\end{figure}

If a state $\S'$ is obtained from a state $\S$ by a counterclockwise transposition at a segment at a segment $j$, then we denote the \df{weight ratio} between $\S'$ and $\S$ by $w(j)$. Thus, \[ w(j) = \frac{w(\S')}{w(\S)}.\]
Note that $w(j)$ depends only on the segment $j$ and not on the state $\S$ and $\S'$.
The possible values for $w(j)$ and its specialization at $W = t^{\frac{1}{2}}$, $B = t^{-\frac{1}{2}}$ are shown in Figure \ref{fig::weightsegment}.

%%%%%%%%%%%%%%%%%%%%%%%%%%%%%%%%%%%%
%% SECTION 
%%
%%%%%%%%%%%%%%%%%%%%%%%%%%%%%%%%%%
\section{The Jacobian algebra of a link diagram}\label{sect B}
Let $K$ be a reduced diagram of an oriented prime  link without curls. Denote by $n$ the number of crossings and label the segments $1,2,\dots,2n$  as in section \ref{sect 1}. 
We shall use the notation $K_0$ for the set of crossing points, $K_1$ for the set of segments, and $K_2$ for the set of regions (including the unbounded region) of $\mathbb{R}^2\setminus K$. 

We are going to construct a quiver with potential and consider its Jacobian algebra.

We define the quiver $Q$  as follows. The set of vertices $Q_0$ is the set of segments of $K$. Thus $Q_0=K_1$. The set of arrows $Q_1$ is the set of arrows of $K$ introduced in section~\ref{sect 1}, more precisely, there is an arrow $i\to j$ in $Q$ if and only if 
\begin{itemize}
\item the segments $i $ and $j$ of $Q$ meet at a crossing point $p$;
\item the segments $i$ and $j$ bound the same region $R$;
\item going clockwise around $p$, we cross $i$ then $R$ then $j$. 
\end{itemize}

\begin{figure}
\begin{center}
\[ \xymatrix@C50pt{3 \ar[rr]_{\zg_3} \ar[d]_{\za_3} && 6\ar@{<-}@<2pt>[d]^{\zd_4}\ar@<-2pt>[d]_{\zg_4} \ar@/^80pt/[ldddd]_{\zd_1} \\
8\ar[r]_{\zb_1}\ar@{<-}@<2pt>[d]^{\zb_4}\ar@<-2pt>[d]_{\za_4}&5\ar[rd]_{\zb_2}\ar[lu]_{\zg_2} & 2\ar[l]_{\zg_1} \\
4\ar[rdd]_{\za_1}&&7 \ar[ll]_{\zb_3}\ar[u]_{\zd_3} \\ \\
&1\ar@/^80pt/[luuuu]_{\za_2} \ar[ruu]_{\zd_2}
}
\]
\caption{The quiver of the figure-eight knot of Example~\ref{ex::kstateslattice}.}
\label{fig:quivereight}
\end{center}
\end{figure}
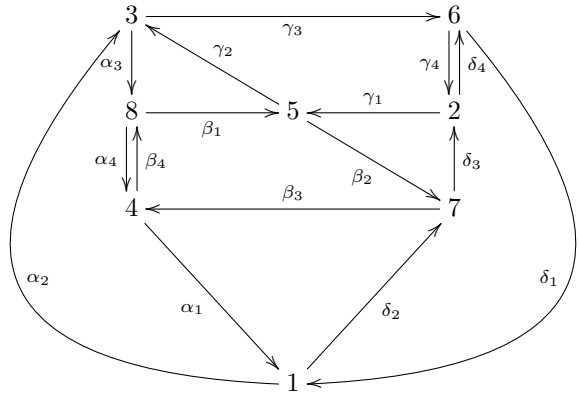

For example the quiver of the figure-eight knot in Example~\ref{ex::kstateslattice} is shown in Figure~\ref{fig:quivereight}. 

The planar  link diagram $K$ induces a planar embedding of $Q$, and since $K$ has no curls, $Q$ has no loops. On the other hand, $Q$ may have 2-cycles, see however section~\ref{sect 2-cycles}. 

The quiver $Q$ has the following two types of chordless cycles.  For each crossing point $p\in K_0$, we obtain a 4-cycle $\omega_P$ and for each region $R$ bounded by $r$ segments, we obtain an $r$-cycle $\omega_R$. 
Each arrow $(p,R)$ lies in exactly two chordless cycles $\omega_p$ and $\omega_R$. 

We define a potential $W$ as
\[ W=\sum_{p\in K_0} \omega_p -\sum_{R\in K_2} \omega_R.\]
In the example in Figure~\ref{fig:quivereight}, the potential is 
\[
\begin{array}
 {rcl}
 W&=&\za_1\za_2\za_3\za_4 +\zb_1\zb_2\zb_3\zb_4+\zg_1\zg_2\zg_3\zg_4+\zd_1\zd_2\zd_3\zd_4
\\
&&-\za_1\zd_2\zb_3-\za_2\zg_3\zd_1-\za_3\zb_1\zg_2-\zb_2\zd_3\zg_1-\za_4\zb_4-\zg_4\zd_4
\end{array}
\]
where the first row consists of the four 4-cycles of the four crossing points in $K$ and the second row consists of the four 3-cycles and two 2-cycles of the six regions in $K$.
\begin{definition}
 The algebra $B=\textup{Jac}(Q,W)$ is called the (completed) \emph{Jacobian algebra of the link diagram $K$}. 
\end{definition}
\begin{remark}
(a) The quiver $Q$ and  the algebra $B$ are not invariants of the link. For example, the second Reidemeister move will change the number of vertices in $Q$.

(b) The quiver does not see the difference between an overpass and an underpass in $K$.
\end{remark}
\subsection{Removal of 2-cycles}\label{sect 2-cycles}
 Each bigon in the  link diagram gives rise to a 2-cycle in the quiver. We can replace the quiver with potential $(Q,W)$ by a quiver with potential $(Q',W')$ without 2-cycles as follows. The quiver $Q'$ is obtained from $Q$ by removing all 2-cycles. The  potential $W'$ is obtained from $W$ as follows.
For every bigon $R$, given by two segments $i,j$ that cross each other in two crossing points $p_1,p_2$, we replace $\omega_{p_1}+\omega_{p_2}-\omega_R$ by the 6-cycle $(\partial_{(p_1,R)} \omega_{p_1})(\partial_{(p_2,R)} \omega_{p_2})$ obtained by joining the two 4-cycles $\omega_{p_1}$ and $\omega_{p_2}$.
This identification on all 2-cycles induces an isomorphism of algebras 
\[B=\textup{Jac}(Q,W) \cong 
\textup{Jac}(Q',W').\] 
This realization of the algebra $B$ by a quiver without loops and 2-cycles will be important when we describe the connection to cluster algebras.

In our running example, there are two bigons formed by the pairs of segments $(4,8)$ and $(2,6)$ in Example~\ref{ex::kstateslattice} and these give rise to two 2-cycles $\za_4\zb_4$ and $\zg_4\zd_4$ in the quiver in Figure~\ref{fig:quivereight}. The above reduction produces the potential
\[
\begin{array}
 {rcl}
 W'&=&\za_1\za_2\za_3\zb_1\zb_2\zb_3+\zg_1\zg_2\zg_3\zd_1\zd_2\zd_3
\\
&&-\za_1\zd_2\zb_3-\za_2\zg_3\zd_1-\za_3\zb_1\zg_2-\zb_2\zd_3\zg_1.
\end{array}
\]

Note that  $\partial_{\za_4} W = -\zb_4 + \za_1 \za_2 \za_3 $ and thus, in the Jacobian algebra $\textup{Jac}(Q,W)$, the arrow $\zb_4$ is equal to the path of length three $ \za_1 \za_2 \za_3$. Similarly  $\za_4$ is equal to $ \zb_1\zb_2\zb_3$.
%%%%%%%%%%%%%%%%%%
%% SECTION
%%%%%%%%%%%%%%%%%%
\section{The link diagram module $T$}\label{sect knot module}
Let $K$ be a curl free diagram of a prime  link and let $B=\textup{Jac}(Q,W)$ be its Jacobian algebra. In this section, we associate a $B$-module $T=\oplus_{i\in K_1} T(i)$ to $K$, where each $T(i)$ is an indecomposable summand. 
The $T(i)$ are constructed explicitly as representations of the quiver $Q$. 

\subsection{A partition of $K_1$}\label{sect partition}
Let $i$ be a fixed segment of $K$. We shall define a partition of the set of all segments $K_1=\sqcup_{d\ge 0} K(d)$ and use it later to define the representation $T(i)$. The sets $K(d)$ depend on the choice of the segment $i$, but, in the interest of simplicity, our notation does not reflect this dependency. This should not create confusion, since $i$ is fixed here. The construction is recursive and the case $d=0$ is slightly different from the cases $d>0$.  But first let us run through the construction in the following example.
\subsubsection{An example}
\begin{figure}
\begin{center}
\scalebox{.8}{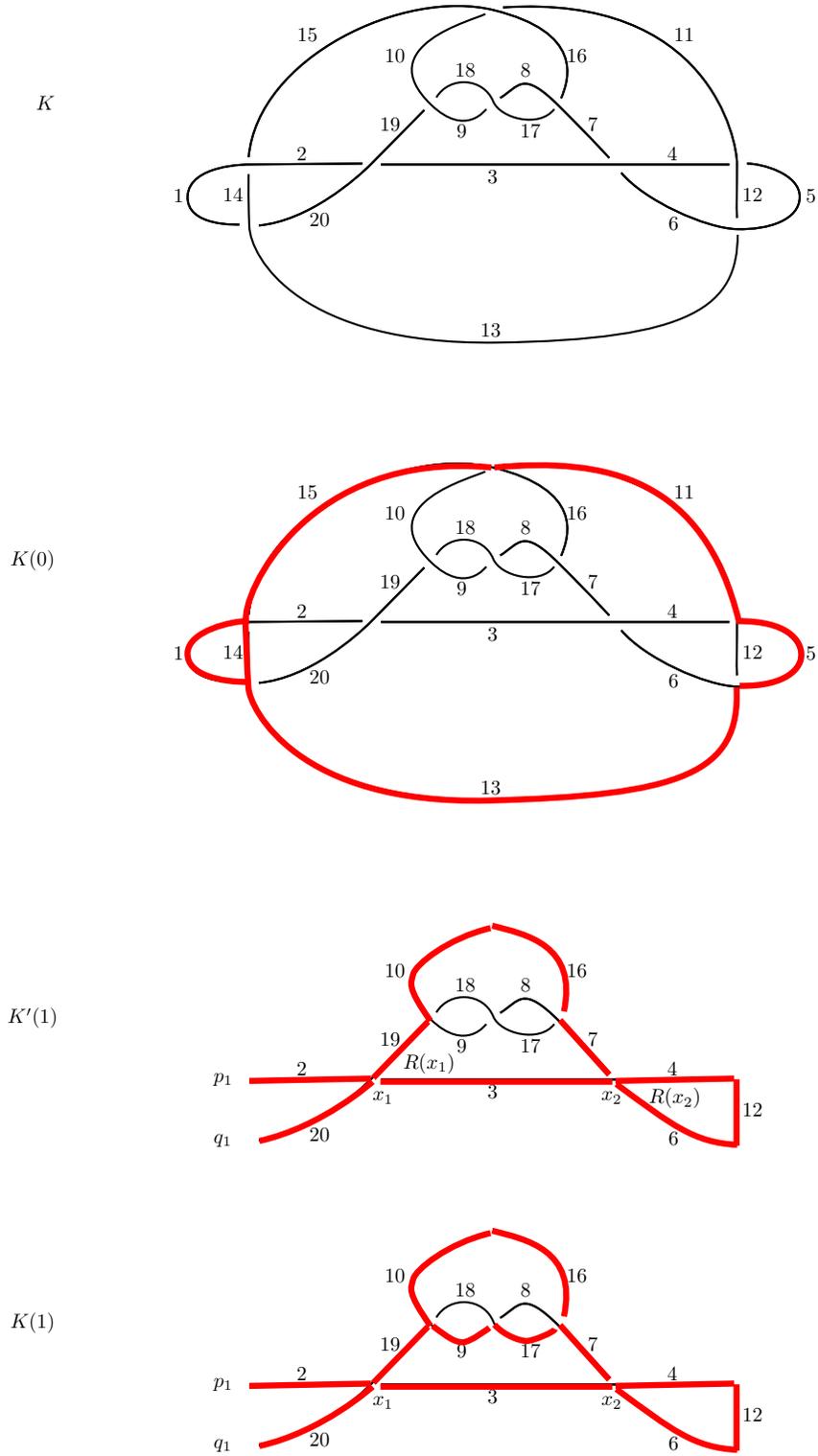}
\caption{An example of the construction of the partition of the segments of the knot into disjoint subsets $K(d)$. The quiver of this diagram is shown in Figure~\ref{fig 515}.}
\label{fig knot 10}
\end{center}
\end{figure}

Consider the knot diagram $K$ illustrated in the top picture in Figure~\ref{fig knot 10}. This is the knot $10_{66}$ in the Rolfsen table.  We choose the segment $i=1$. The set $K(0)$ is the set of all edges that share a region with the segment $i=1$, including $i$ itself. This set is shown in red in the second picture in the figure. Thus we have $K(0)=\{1, 15,11,5,13, 14\}$. We think of this set as a union of two paths both starting and ending with the segment 1.
The first path $w_{L,0}$ starts on segment 1 in the direction given by the orientation of the knot and turns left at each crossing point until it comes back to $1$. Thus $w_{L,0}= 1,15,11,5,13,1$. The second path
$w_{R,0}$ also starts on segment 1 in the same direction, but it turns right at each crossing point. Thus $w_{R,0}=1,14,1$. 

The set $K(1)$ is constructed in two steps. First, we remove the set $K(0)$ from $K$, and then we define the set $K'(1)$ as the set of  all segments  that are incident to the unbounded region of $K\setminus K(0)$.  This set is shown in red in  the third picture in Figure~\ref{fig knot 10}. Note that there are precisely two crossing points $p_{1}$ and $q_{1}$ that are incident to exactly one segment of $K'(1)$. Again, we can think of this set as the union of two paths, but this time they start at $p_{1}$ and end at $q_{1}$.
The first path $w_{L,1}$ makes a left turn at every crossing point. Thus $w_{L,1}=2,19,10,16,7,4,12,6,3,20.$ The second path $w_{R,1}$ makes a right turn at every crossing point. Thus $w_{R,1}=2,20.$

In our example there are two crossing points $x_1$ and $x_2$ that are of degree 4 in $K'(1)$. In this situation, the set $K(1)$ is strictly larger than $K'(1)$. It is shown in the last picture in the figure and is defined as follows. 
The path $w_{L,1}$ goes through each of the points $x_1,x_2$ exactly twice. Let $D(x_i)$ denote the domain in the plane bounded by the subpath of $w_{L,1}$ from $x_i$ to $x_i$. Thus $D(x_1) $ is bounded by the path $19,10,16,7,4,12,6,3$ and $D(x_2) $ is bounded by the path $4,12,6$. The domain $D(x_2)$ actually consist of a single region of the  diagram. On the other hand, the domain $D(x_1)$ contains 5 regions of $K$. We let $R(x_i)$ be the unique region of $K$ inside $D(x_i)$ that is incident to $x_i$. Then $R(x_2)=D(x_2)$ and $R(x_1)$ is the region bounded by the segments $19,9,17,7,3$. 

Then $K(1)$ is defined as the union of $K'(1)$ with the segments of the regions $R(x_i)$. Thus we need to add the segments $9$ and $17$ to  our set. We are now done with the case $d=1$. 

The set $K(2)$ is again defined in two steps, but the second step will be trivial.  First  let $K'(2)$ be the set of all segments that are incident to the unbounded region of $K\setminus (K(0)\cup K(1))$.
Thus $K'(2)=\{18,8\}$. There are no crossing points of degree 4 in this set, and therefore we have $K(2)=K'(2)$. 

\subsubsection{The general case for $d=0$}
For a general link, define
\[K'(0)=K(0)=\{j\in K_1\mid \textup{$j$ and $i$ bound the same region of $K$}\}\cup\{i\},\]
and let $\overline{K'(0)} $ be the closure of $K'(0)$; here closure means that the set also contains the endpoints of the segments.

We can describe $K(0)$ as the union of two paths given by the boundaries of the two regions incident to $i$. We describe these paths below in a way that will generalize to an iteration of this procedure to $d>0$. The set $K(0)$ can be described as the union of the segments along two paths
\[ 
\begin{array}{c}
\xymatrix@C40pt{w_{L,0}\colon p_{0}\ar@{-}[r]_(.6){w_{p,0}} &p'_{0} \ar@{-}[r]_{w'_{L,0}} & q'_{0}\ar@{-}[r]_{w_{q,0}} &q_{0}}
\\
\xymatrix@C40pt{w_{R,d}\colon p_{0}\ar@{-}[r]_(.6){w_{p,0}} &p'_{0} \ar@{-}[r]_{w'_{R,0}} & q'_{0}\ar@{-}[r]_{w_{q,0}} &q_{0}}
\end{array}
\] 
as follows.
Let  $p_{0}$ and $q_{0}$ be the endpoints of the segment $i$. Since $K$ has no curls, we have $p_{0}\ne q_{0}$. Define $p_{0}'=q_{0} $ and $q_{0}'=p_{0}$ and let $w_{p,0}=w_{q,0}$ be the segment $i$. At $p_{0}'$ and at every subsequent crossing point, the path $w_{L,0}$ turns left, and therefore $w_{L,0}$ is the boundary of the region to the left of the segment $i$ from $p_{0}$ to $q_{0}$. Similarly $w_{R,0}$ is the boundary of the region to the right of the segment $i$. These are exactly the segments in $K'(0)$. 

Note that the two points $p_{0},q_{0}$ are of degree 3 in $\overline{K(0)}$ and all other crossing points in $K$ have degree 0 or 2  in $\overline{K(0)}$.

\begin{lemma}\label{lem d=0}
The two subpaths $w'_{L,0}$ and $w'_{R,0}$ do not share a crossing point besides $p'_{0}$ and $q'_{0}$. 
\end{lemma}
\begin{proof}
 Let $R_1,R_2\in K_2$ denote the two regions at $i$. Suppose there exists $p\in K_0$ such that $p$ is not an endpoint of $i$ and $p\in R_1\cap R_2$. Then we can draw a closed curve $\zg$ from $p$ to $p$ that runs through $R_1,R_2$, crosses $i$ once, and does not cross any other segment of $K$. We consider two cases, depending on the local configuration of the four segments at $p$ relative to $\zg$, see Figure~\ref{fig gamma}.  
 \begin{figure}
\begin{center}
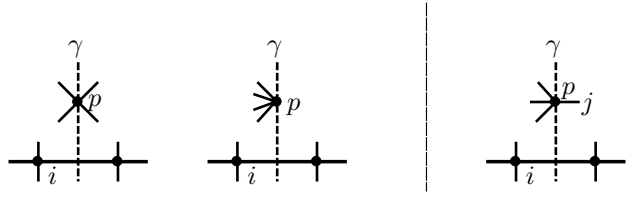
\caption{Proof of Lemma \ref{lem d=0}.
At the point $p$, there is an even number of edges on either side of the curve $\zg$ in the two pictures on the left, and an odd number  in the picture on the right.}
\label{fig gamma}
\end{center}
\end{figure}

Suppose first that there is an even number of these four segments on either side of $\zg$. This case is illustrated in the two pictures on the left of the figure. Note that at either endpoint of the segment $i$, there are three loose segments, so that in total there is an odd number of loose segments on either side of $\gamma$, and (using the Jordan curve theorem) it is thus impossible to pair them up without crossing $\zg$ in order to form a  link.

Therefore, out of the four segments at $p$, there must be an odd number on either side of $\zg$. This case is illustrated in the right picture of the figure. Then, on one side, there is only one segment; call it $j$. Moving $\zg$ slightly away from $p$ toward the segment $j$, we obtain a closed curve $\zg'$ that crosses only two segments $i$ and $j$. This shows that $K$ is the connected sum of two  links and thus not prime, a contradiction.
\end{proof}

\subsubsection{The general case for $d\ge 1$}
We define $K(d)$ recursively. 

\begin{definition}\label{def K'}
Assume $K(e)$ is defined for all $e<d$. 
Let $K'(d)$ be the set of all segments $j$ in $K_1\setminus \left(\cup_{e<d}\, K(e)\right)$ for which there exists a segment $k\in K(d-1)$  such that $j$ and $k$ bound the same region of $K$.
Let $\overline{K'(d)}$ be the closure of $K'(d)$.

 \end{definition}

\begin{definition}\label{def R(x)}
 (a) A crossing point $p\in K_0$ is called  an \emph{external point in $K'(d)$} if exactly one of its incident segments lies in $K'(d)$. A segment $j\in K_1$ is called \emph{external in $K'(d)$} if $j\in K'(d)$ and exactly one the endpoints of $j$ is an external point in $K'(d)$.

(b) A crossing point $p$ is called an \emph{internal point} in $K'(d)$ if all four segments at $p$ lie in $K'(d)$. Note that for $d=0$ there are no internal crossing points since our diagram $K$ has no curls. 
If there exists a non-constant path $w$ starting and ending at $x$ that uses only segments of $K'(d)$, we let $D(x)$ be the bounded domain enclosed by $w$ in the plane. 
Then $D(x)$ is a union of regions of $\overline{K'(d)}$. 
 Let $R(x)\in K_2$ be the unique region in $D(x)$ that contains $x$.  
\end{definition}
 
\begin{figure}
\begin{center}
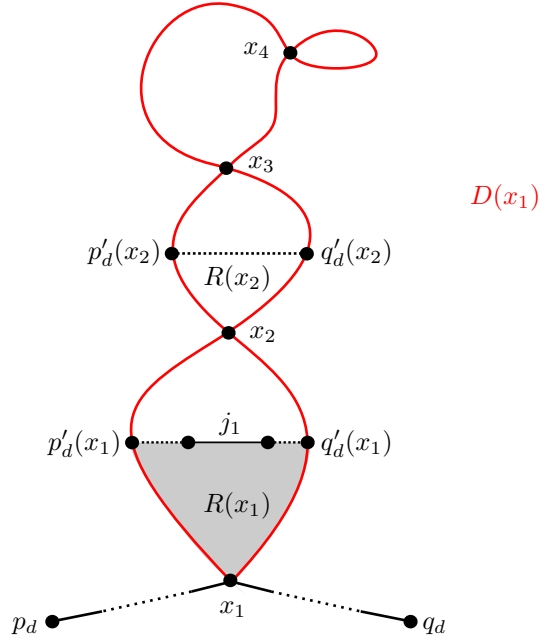
\caption{The domain $D(x_1)$ and the region $R(x_1)$ associated to an interior point $x_1\in K(d)$. The domain $D(x_1)$ is bounded by the red curve and the region $R(x_1)$ is shaded. It is bounded above by the black segments, one of which is labeled $j_1$. The domain $D(x_2)$ for the second interior point $x_2$ is only the part of $D(x_1)$ that lies above the point $x_2$.}
\label{fig D(x)}
\end{center}
\end{figure}

 An example is shown in Figure~\ref{fig D(x)}. In that figure, there are four interior points $x_1,\dots,x_4$. The domain $D(x_1)$ is bounded by the red subcurve $w$ and the region $R(x_1)$ is shaded. 

 Now let $j$ be a segment of $K$. 
 We define $\ze_d(j)\in\{0,1\}$ for $d\ge 1$ as follows. 
 \begin{equation}\label{eq epsilon}
\ze_d(j) = \left\{
\begin{array}
 {ll}
 1 & \textup{if there exists an internal point $x$ in $K'(d)$ such that $j$ lies} \\ &\textup{in the interior of $D(x)$ and the region $R(x)\in K_2$ contains the} \\ &\textup{segment $j$ and the point $x$;}
\\
0 & \textup{otherwise.}
\end{array}
\right.
\end{equation}

\begin{definition}
For $d\ge 1$, let
 \[K(d) =K'(d)\cup\{j\in K_1\mid \ze_d(j)=1\}.\]
 \end{definition}

For example, the segment $j_1$ in Figure~\ref{fig D(x)} satisfies the first condition for $x=x_1$. Thus $\ze_d(j_1)=1$ in this case.
The set $K(d)$ contains every segment of the red curve and every segment of the black curves bounding $R(x_1)$ and $R(x_2)$.

\begin{lemma}
 \label{lemma 5} Let $d\ge1$.

(a) Each connected component $C$ of $\overline{K'(d)}$ is either a single path $w$ from $ p_{d}  $ to $q_{d}$ or
the union of two paths
\[ 
\begin{array}{c}
\xymatrix@C40pt{w_{L,d}\colon p_{d}\ar@{-}[r]_(.6){w_{p,d}} &p'_{d} \ar@{-}[r]_{w'_{L,d}} & q'_{d}\ar@{-}[r]_{w_{q,d}} &q_{d}}
\\
\xymatrix@C40pt{w_{R,d}\colon p_{d}\ar@{-}[r]_(.6){w_{p,d}} &p'_{d} \ar@{-}[r]_{w'_{R,d}} & q'_{d}\ar@{-}[r]_{w_{q,d}} &q_{d}}
\end{array}
\]
from $ p_{d}  $ to $q_{d}$, where  
 $ p_{d}  $ and  $q_{d}$ are external points in $K'(d)$,
the initial and terminal subpaths $w_{p,d}$ and $w_{q,d}$ are the same in both paths, $\deg p_{d}'=\deg q_{d}'=3$ in $C$, $p_{d}'\ne q_{d}'$ and $w'_{L,d}$ (respectively $w'_{R,d}$) is obtained by turning left (respectively right) at every crossing point in $K'(d)$.

If no such points $p'_{d},q'_{d}$ exist then $ w_{L,d}=w_{R,d}$ and $K'(d)$ is a single path from $p_{d}$ to $q_{d}$.

(b) All other crossing points, besides $p_{d},p
'_{d},q_{d},q'_{d}$, have degree 0, 2 or 4 in $C$. 
In particular, $C$ has exactly two external points $p_{d} $ and $q_{d} $, and moreover, $p_{d} ,q_{d} \in \overline{K'(d)} \cap \overline{K'(d-1)}$ and $p_{d},q_{d}\notin \overline{K'(e)}$, with $e\ne d, d-1$.

(c)  The paths $w'_{L,d}$ and $w'_{R,d}$ do not share a crossing point besides $p_{d}'$ and $q_{d}'$.
\end{lemma}
\begin{proof}
(a) Suppose first $d=1$. The external points are $p_{1}=p'_{0}$ and
 $q_{1}=q'_{0}$, and there are no other external points in $K'(1)$. 
  If there are points of degree 3 in $\overline{K'(1)}$, we let $p_{1}'$ be the point of degree 3 closest to $p_{1}$, and let $q_{1}'$ be the point of degree 3 closest to $q_{1}$. Note that every connected component has an even number of points of degree 3, because of parity, and that there are at most two  because there are only two external points.
  
   Let $w_{p,1}$ be the unique path from $p_{1}$ to $p'_{1}$ in $K'(1)$ that
\begin{equation}\label{condition path}
\begin{array}{l}
\bullet\   \textup{does not use an edge twice;}\\
\bullet \ \textup{is of maximal length;}\\
\bullet \ \textup{turns left or right at every point of degree 4 in $\overline{K'(d)}$.}
\end{array}
\end{equation}
Such a path exists and is unique by the following argument. The starting point $p_{1}$ is of degree 1 in $C$, so the first step is uniquely determined. At every point of degree 2, the incoming edge leaves only one choice for the outgoing edge. At a point of degree 4, there are a priori two possibilities, turn left or turn right, but only one of these will produce a path of maximal length. 
%Note that if there is an internal point $x$ in $C$ then the boundary of the domain $D(x)$ lies entirely in one of the two paths $w_{L,1}$ or $w_{R,1}$.
Similarly, let $w_{q,1}$ be the unique path from $q'_{1}$ to $q_{1}$ in $K'(1)$ that
respects conditions (\ref{condition path}).

Then the paths $w_{L,1},w_{R,1}$ form the boundary of the regions on the left and right of the path $(w_{q,1}\, i\, w_{p,1})$ in $(K\setminus K(0))\cup\{i\}$. These are exactly the segments in $K'(1)$.
This completes the proof of (a) for $d=1$.

(b)  The degree formulas follow directly from (a). In particular $p_{1}$ and $q_{1}$ are the only external points in $C$. Furthermore, three of the segments at $p_{1}$ lie in $K(0)$, and the remaining segment, which is the first segment of $w_{p,1}$, lies in $K'(1)$. 
Thus $p_{0}\in \overline{K'(1)}\cap \overline{K'(0)}$ and $p_{0}\notin K'(e)$, with $e\ne 0,1$. The proof for $q_{1}$ is similar. 
This also implies that all other points have degree 0, 2 or 4 in $\overline{K'(1)}.$

(c) Suppose $w'_{L,1}$ and $w'_{R,1}$ share a point $x\ne p_{1}',q_{1}'$, see Figure~\ref{fig lem 5c}.
\begin{figure}
\begin{center}
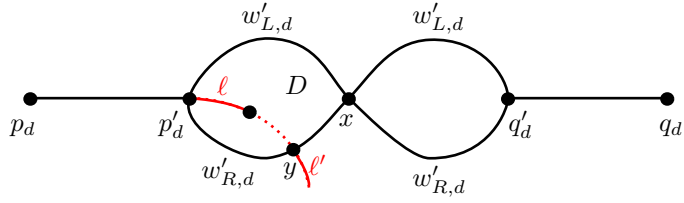
\caption{Proof of Lemma \ref{lemma 5} part (c).}
\label{fig lem 5c}
\end{center}
\end{figure}
Let $D$ denote the domain in the plane bounded by the segments of the paths $w_{L,1}$ and $w_{R,1}$ from $p'_{1}$ to $x$.
Let $\ell$ denote the segment at $p_{1}'$ that does not belong to $w_{L,1}$ or $w_{R,1}$. Since $w_{L,1}$ turns left at $p_{1}'$ and $w_{R,1}$ turns right, the segment $\ell$ lies in $D$.
Following the  link, starting at $p'_{1}$ in direction of $\ell$, we must reach a point $y$ where we leave $D$. Let $\ell'$ denote the segment outside $D$ right after $y$. In $K$, this segment $\ell'$ bounds the same region as a segment on the path $w'_{R,1}$, and therefore $\ell'$ must lie in $K(0)\cup K(1)$.
However, $\ell' $ cannot lie in $K(0)$, because $K(0)$ has no external segments.  On the other hand, $\ell'$ cannot be in $K'(1)$, because (a) implies that every segment of $K'(1)$ lies on one of the two paths $w_{L,1}, w_{R,1}$.
This is a contradiction, and thus the two paths $w'_{L,1},w'_{R,1}$ cannot have the point $x$ in common. 
This completes the proof for $d=1$.

For $d>1$, the proof is similar, with the additional feature that now we may have components that arise from internal points in $K'(d-1)$. Indeed, every internal point $x$ of $K'(d-1)$, such that there is a segment $j$ bounding the region $R(x)$ with $\ze_{d-1}(j)=1$, determines two crossing points $p'_{d}(x)$ and $q'_{d}(x)$ as the unique points in $R(x) \setminus \{x\}$ that lie on the boundary of $D(x)$, and such that the path $w=w_{L,d-1}$ or $w_{R,d-1}$ is of the form
 \[\xymatrix{w\colon p_{d}\ar@{-}[r] & x\ar@{-}[r]& p_{d}'(x) \ar@{-}[r]
& q_{d}'(x) \ar@{-}[r]
 &x \ar@{-}[r]
 &q_{d}},
 \]
 see Figure~\ref{fig D(x)}.
Therefore every internal point $x$ of $K'(d-1)$ gives rise to a connected component of $K(d)$ in which the  points $p_{d}=p'_{d-1}(x)$ and $q_{d}=q'_{d-1}(x)$ are the unique external points. For all other components, we have $p_{d}=p'_{d-1}$ and $q_{d}=q'_{d-1}$.
 Note that  $p_{d}$ and $q_{d}$ have degree 3 in $\overline{K(d-1)}$ and thus degree 1 in $\overline{K'(d)}$. Hence they are external points.

The rest of the proof of (a) is analogous to the case $d=1$. The point $p_{d}'$ is the point of degree 3 closest to $p_{d}$, and 
$q_{d}'$ is the point of degree 3 closest to $q_{d}$.
The paths $w_{p,d}$ and $ w_{q,d}$ are the unique paths from $p_{d}$ to $p_{d}'$ and from $q_{d}$ to $q_{d}'$ that satisfy the conditions (\ref{condition path}). Moreover the paths $w_{L,d}, w_{R,d}$ form the boundaries of the regions to the left and right of the path
\[ w_{q,d}\dots w_{q,2} w_{q,1} \, i\, w_{p,1} w_{p,2} \dots w_{p,d}\]
in 
$(K\setminus   \cup_{e=0}^{d-1} K(e) ) \cup \{ w_{q,d-1},\dots w_{q,1}, i, w_{p,1}, \dots, w_{p,d-1}\}$.
Thus $w_{L,d}$ and $w_{R,d}$ consist of the edges of $K'(d)$.

 The proof of (b) is analogous to the case $d=1$ except that the point $p_{d}$ may be equal to a point $p'_{d-1}(x)$ for some interior point $x$ in $K'(d)$.  

In the proof of (c), the only difference to the proof in case $d=1$ is that now the segment $\ell'$ in Figure~\ref{fig lem 5c} cannot lie in $ K(d-1)$, because otherwise $y$ would be an external point of $K(d-1)$ that lies in $\overline{K(d)}$, a contradiction to (b). 
%%
%%%
%%%
%%%
%(c) Suppose $w'_{L,d}$ and $w'_{R,d}$ share a point $x\ne p_{d}',q_{d}'$, see Figure~\ref{fig lem 5c}.
%\begin{figure}
%\begin{center}
%\input{figlem5c.pdf_tex}
%\caption{Proof of Lemma \ref{lemma 5} part (c).}
%\label{fig lem 5c}
%\end{center}
%\end{figure}
%Let $D$ denote the domain in the plane bounded by the segments of the paths $w_{L,d}$ and $w_{R,d}$ from $p'_{d}$ to $x$.
%Let $\ell$ denote the segment at $p_{d}'$ that does not belong to $w_{L,d}$ or $w_{R,d}$. Since $w_{L,d}$ turns left at $p_{d}'$ and $w_{R,d}$ turns right, the segment $\ell$ lies in $D$.
%Following the  link, starting at $p'_{d}$ in direction of $\ell$, we must reach a point $y$ where we leave $D$. Let $\ell'$ denote the segment outside $D$ right after $y$. In $K$, this segment $\ell'$ bounds the same region as a segment on the path $w'_{R,d}$, and therefore $\ell'$ must lie in $K(d-1)\cup K(d)$.
%If $\ell'\in K(d-1)$ then $y$ is an external point of $K(d-1)$ that lies in $\overline{K(d)}$, a contradiction to (b). On the other hand, $\ell'$ cannot be in $K'(d)$, because (a) implies that every segment of $K'(d)$ lies on one of the two paths $w_{L,d}, w_{R,d}$.
%This is a contradiction, and thus the two paths $w'_{L,d},w'_{R,d}$ cannot have the point $x$ in common. 
\end{proof}

\begin{remark}
 Each interior point $x$ of $K'(d-1)$ whose region $R(x)$ contains a segment $j$ with $\ze_{d-1}(j)=1$ gives rise to a connected component of $K'(d)$.
\end{remark}

\begin{proposition}\label{prop partition}
For every segment $i\in K_1$, we have constructed a partition \[K_1=\sqcup_{d\ge 0} K(d).\]
\end{proposition}
\begin{proof}
 This follows directly from the construction.
\end{proof}

We are now ready to define the dimension vector of the representation $T(i)$.
\begin{definition} \label{def dim}
Let $K_1=\sqcup_{d\ge 0} K(d) $ be the partition with respect to a segment $i\in K_1$. For every  segment $j\in K_1$, we define 
$d(i)_j = d $ if $j\in K(d)$.
\end{definition}

In the example of Figure \ref{fig knot 10}, we have $d(1)_j =1$ if $j=2,3,4,6,7,9, 10,12,16,17, 19 , 20$; $d(1)_j=2$ for $j=8,18$, and $d(i)_j=0$ for all other $j$.

Our next result says that the dimension difference at adjacent vertices is at most one.
\begin{proposition} \label{prop dim}
Let $i,j,k\in K_1$. 
 If there is an arrow $j\to k \in Q$ then 
 \[|d(i)_j-d(i)_k|\le 1.\]
\end{proposition}
\begin{proof}
Let $d=d(i)_j$. Thus $j\in K(d)$. The existence of the arrow $j\to k$ implies that $j$ and $k$ bound the same region in $K$. If $k\notin \cup_{e\le d} K(e)$ then Definition~\ref{def K'} implies that $k\in K'(d+1)$, thus $k\in K(d+1) $ and $d(i)_j-d(i)_k=-1$. 
If $k\in K(d)\cup K(d-1)$ then $d(i)_k\in \{d,d-1\}$ and there is nothing to show.
Finally suppose $k\in K(e) $, with $e\le d-2$. Then Definition~\ref{def K'}  implies that $j$ lies in $K'(e+1)$ unless it already lies in $\cup_{e'<e} K(e')$. In both cases, we  have $d(i)_j\le e+1\le d-1$, a contradiction.
\end{proof}

\subsection{Properties of $K(d)$}

It will be convenient to use the following terminology.
Given two segments $i,j\in K_1$, a curve in $\mathbb{R}^2$ is called a \emph{dimension curve} from  $j$ to  $i$ if it starts at a point on segment $j$, ends at a point on segment $i$ and does not go through a crossing point of $K$. 

Let $\dimprime(i,j)$ be the minimal number of crossings between the segments of $K$ and a dimension curve from segment $j$ to segment $i$. 
We call a segment $j\in K_1$ an \emph{interior segment of} $K(d)$ if $K(d)$ contains all the segments on the boundary of the two regions incident to $j$ in $K$.

Note that the segment $i$ is an interior segment of $K(0)$. The following lemma says that there are no other interior segments.
\begin{lemma}
 \label{lemma 6} If $d\ge 1$ then
 $K(d)$ has no interior segments.
\end{lemma}
\begin{proof}
Suppose a segment $j$ belongs to two regions $R_1,R_2$ in $K$ and each segment in $R_1\cup R_2$ lies in $K(d)$, see Figure~\ref{fig lem 6}. 
\begin{figure}
\begin{center}
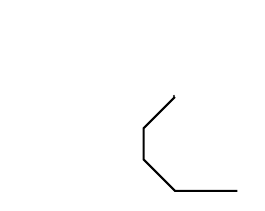
\caption{Proof of Lemma \ref{lemma 6}.}
\label{fig lem 6}
\end{center}
\end{figure}
Since $d\ge 1$, the dimension curve of $j$ must cross a segment $k$ of $R_1\cup R_2$, so it has one more crossing than the dimension curve of the segment $k$. Since $j,k\in K(d)$ this means that $\ze_d(j)=1$. Thus there exists an internal point $x$ in $K(d)$ satisfying the condition (\ref{eq epsilon}). In particular, 
one of the two regions at $j$, say $R_1$, contains $x$ and $j$. Thus $R_1$ is the region $R(x)$ of condition~(\ref{eq epsilon}).
 The other region of $K$ at $j$ is the region $R_2$ and both lie entirely in $K(d)$. 
Since $j\notin K'(d)$, the segments of $R_1\cup R_2$ that do lie in $K'(d)$ all lie in the same region in $K'(d)$.  
In particular, no segment $k$ of $R_1\cup R_2$ shares a region with another internal point $x'\ne x$ such that $k$ lies in $D(x')$. Therefore each segment $k$ of $R_2\setminus\{j\}$ satisfies $\ze_d(k)=0$. Thus  $k\in K'(d)$, which implies $k\in D(x)$.

Now consider an endpoint $y$ of $j$. Three of the segments incident  to $y$ lie in  $R_1\cup R_2$ and the fourth segment $\ell$ doesn't, see Figure~\ref{fig lem 6}. Therefore $\ell$ either lies on the boundary of $D(x)$ and hence $\dimprime(i,\ell)=d$, or $\ell$ lies outside $D(x)$ and hence $\dimprime(i,\ell)=d-1$. The latter case is impossible by Lemma~\ref{lemma 5}(b), because $\ell$ would be an external segment of $K'(d-1)$, but $\ell$ has an endpoint in $\overline{K'(d)}$ (and not in $\overline{K'(d-2)}$\ ). 

Thus the segment $\ell$ lies on the boundary of $D(x)$. Since $\ell$ is not in $R_1$, there exists a point $y'\in K_0$ such that two of its incident segments lie in $R_1$ and one, call it $h$, lies on the boundary of $D(x)$ between $x$ and $y$ as in the figure. 
Denote by $\ell'$ the fourth segment at $y'$. Note that it must lie inside $D(x)$, because otherwise $y'$ would be an external point of ${K'(d-1)}$ that does not lie in $\overline{K'(d-2)}$, contradicting Lemma~\ref{lemma 5}(b).
Since $R_1$ is a region in $K$, the segment $\ell$ must lie in the connected component $C$ of $D(x)\setminus R_1$ that contains $y'$. Following the  link starting at $y'$ in direction $\ell'$, we must reach a point $y''$, where we leave the component $C$. Let $\ell''$ be the segment outside $C$ right after $y''$. Then $\ell''$ is an external segment of $K'(d-1)$ with external point $y''$ not in $\overline{K(d-2)}$, again a contradiction to Lemma~\ref{lemma 5}(b).
\end{proof}

It will be convenient to consider the following dual graph. For an illustration, see Example~\ref{ex 515}
\begin{definition}
(a) Let $Q(d)$ be the full subquiver of $Q$ on the vertices $j$ such that $j\in K(d)$.

(b) Let $G(d)$ be the graph with vertex set the set of chordless cycles in $Q$ that also lie in $Q(d)$, and two chordless cycles are connected in $G(d)$ by an edge if they share an arrow in $Q(d)$.
\end{definition}

Notice that the graph $G(d)$ has two types of vertices, the crossing point vertices and the region vertices. The first type corresponds to the chordless 4-cycles $\omega_p$, with $p\in K_0$, and the second type corresponds to the chordless cycles $\omega_R$, with $R\in K_2$. 

\begin{corollary}
 \label{cor 7}
 In $G(d)$, the degree of a crossing point vertex is at most 2.
\end{corollary}
\begin{proof}
 If a crossing point vertex $x$ has three adjacent regions in $G(d)$, hence in $K(d)$, then $K(d)$ has an interior segment, contradicting Lemma~\ref{lemma 6}.
\end{proof}

\begin{definition}
 Let $R$ be a region of $K$ such that each segment of $R$ lies in $K(d)$. By Lemma \ref{lemma 5}, one of the two paths $w=w_{L,d}$ or $ w_{R,d}$ encloses the region $R$.  Let $x(R)$ be the first crossing point of the region $R$ on the path $w$.  We call $x(R)$ the \emph{root} of the region $R$.
\end{definition}

Recall that a \emph{leaf} in a graph is a vertex of degree one.
\begin{lemma}
 \label{lemma 8} 
 (a) The mapping $R\mapsto x(R)$ is a bijection between 
the sets of 
region vertices of $G(d)$ and
 crossing point vertices of $G(d)$.

(b) Every connected component of $G(d)$ has a unique vertex that is a leaf and a crossing point vertex.
\end{lemma}
\begin{proof}
 (a) Since the path $w$ starts and ends at vertices outside $R$ it must go through the point $x(R)$ twice, in the sense that it contains all four segments at $x(R)$. Thus $x(R)$ is an internal point of $K'(d)$ and therefore a vertex of $G(d)$. This shows that the mapping is well-defined.
 
 The mapping is injective by definition. Now let $x$ be any crossing point vertex in $G(d)$. Then all four segments at $x$ lie in $K(d)$. This implies that $x$ is an internal point of $K'(d)$, because the endpoints of the segments $j$ with $\ze_d(j)$ have at most degree 3 in $K(d)$. Lemma~\ref{lemma 5} then implies that the four segments at $x$ all lie on one path $w=w_{L,d}$ or $w=w_{R,d}$, and thus $w$ goes through $x$ twice. Then $w $ is of the form 
 \[\xymatrix{w\colon p_{d}\ar@{-}[r]&x\ar@{-}[r]_{w(x)}&x\ar@{-}[r]&q_{d}}\]
 and the subpath $w(x)$ forms the boundary of the domain $D(x)$ in  condition (\ref{eq epsilon}). 
 There is a unique region $R(x)$ of $K$ that lies within $D(x)$ and contains $x$. By definition of $K(d)$, all segments of $R$ lie in $K(d)$. Thus $R$ is a region vertex of $G(d)$ and $x(R)=x$. This shows that the mapping is surjective.
 
 (b) Let $C$ be a connected component of $G(d)$. By Lemma~\ref{lemma 5}, there is a path $w=w_{L,d}$ or $w=w_{R,d}$ that encloses all regions of the region vertices of $C$. Let $x$ be the first point on $w$ that is a crossing point vertex of $C$. Then $x$ is a leaf of $C$.
 
 To show that there is no other crossing point that is a leaf, note first that every crossing point vertex of $G(d)$ is of degree at most 2 in $G(d)$, by Corollary \ref{cor 7}. Now we proceed by induction on the number of region vertices. If there is only one region vertex $R$ in $C$ then $C=\xymatrix{R\ar@{-}[r]&x(R)}$ and we are done. 
 Suppose there is more than one region vertex. Take a leaf $\ell$. If $\ell=R$ is a region vertex, then $C\setminus\{\ell\}$ has $x(R)$ as a leaf and thus $C\setminus\{\ell,x(R)\}$ is connected, and by induction it has no other crossing point vertex that is a leaf than $x$.
On the other hand, if $C$ contains no leaf that is a region vertex, then there are more crossing point vertices than region vertices in $C$, which is impossible by part (a). 
\end{proof}

\begin{remark}
 We don't know if $G(d)$ is a forest.
\end{remark}

\subsection{Definition of the link diagram module $T$}

Let $\kb$ be an algebraically closed field. Let $K$ be an oriented diagram without curls of a prime  link with $n$ crossings. Let $(Q,W)$ be the associated quiver with potential and $B=\textup{Jac(}Q,W)$ its Jacobian algebra. Let $I_d$ denote the identity matrix of rank $d$. We define the \emph{link diagram module} 
\[T=\oplus_{i\in K_1} T(i)\]
of $K$ as follows. 

For each segment $i$ of $K$, the direct summand $T(i)=(T(i)_j,T(i)_\za)_{j\in Q_0,\za\in Q_1}$
is the representation of $Q$ given by 
\[T(i)_j = \kb^{d(i)_j},\] 
for each vertex $j$, where $d(i)_j$ is the dimension defined  in Definition~\ref{def dim};
and for each arrow $\za\colon j\to \ell$, we define the corresponding linear map
\[T(i)_\za =\left\{ 
\begin{array}
 {lll}
  \left(\begin{array}{c|ccc}
       0  & \multicolumn{3}{c}{\multirow{3}{*}{{\scalebox{1.5}{$I_{d-1}$}}}} \\
    \raisebox{2pt}{\vdots} & & &\\
    0 & & & 
  \end{array}\right)
&&\textup{if $d(i)_j=d(i)_\ell +1=d$}
\\
\\
  \left(\begin{array}{cccc}
  \multicolumn{3}{c}{\multirow{3}{*}{\raisebox{0mm}{\scalebox{1.5}{$I_{d-1}$}}}} 
  \\ 
    \\
    \\ \hline
0&\cdots&0 
  \end{array}\right)
&&
\textup{if $d(i)_j+1=d(i)_\ell=d$}\\
\end{array}
\right.
\]
and if $d(i)_j=d(i)_\ell=d$ then

\[T(i)_\za =\left\{ 
\begin{array}
 {cll}
  \left(\begin{array}{c|ccc}
    0  & \multicolumn{3}{c}{\multirow{3}{*}{{\scalebox{1.5}{$I_{d-1}$}}}} \\
    \raisebox{0pt}{\vdots} \\0 & & &\\ \hline
    0 &0 &\cdots & 0
  \end{array}\right)
& &{\genfrac{}{}{0pt}{0}
{\textup{if $\za$ is of the form $(x,R(x))$ with}}
{\textup{ $x$ an internal point of $K'(d)$}}}
\\
\\
\scalebox{1.5}{$I_{d}$}
  &
&\textup{otherwise.}\\
\end{array}
\right.
\]
Because of Proposition \ref{prop dim}, there are no other possibilities for the dimensions and thus $T(i)_\za$ is well-defined.

\begin{example}\label{ex 515}
 The quiver $Q$ of the knot $10_{66} $ in Figure~\ref{fig knot 10} is shown in the top picture and the representation $T(1)$ in the bottom picture of Figure~\ref{fig 515}. The quiver $Q(1)$ is the full subquiver on the vertices 2,3,4,6,7,9,10,12,16,17,19,20. It contains four chordless cycles, the two crossing point cycles $w_{x_1}$, $w_{x_2}$ the two region cycles $w_{R(x_1)}$, $w_{R(x_2)}$, where we use the notation of Figure~\ref{fig knot 10}. Therefore its dual graph is 
 \[\xymatrix{w_{x_1}\ar@{-}[r]&w_{R(x_1)}\ar@{-}[r]&w_{x_2}\ar@{-}[r]&w_{R(x_2)}}.\]

\end{example}

\begin{figure}
\begin{center}
\[\xymatrix@R=19pt@C=19pt{
&& 15 \ar[dddd]\ar[rrrrrr] &&&&&& 11 \ar[ld]\ar[dddddrr] \\
&&& 10\ar[lu]\ar[rd] &&&& 16\ar[dd]\ar[llll]\\
&&&& 18 \ar@<-2pt>[d]\ar[rr]&& 8\ar@<-2pt>[d]\ar[ru]\\
&&& 19\ar[uu]\ar[rrd]&9\ar@<-2pt>[u]\ar[l]&&17\ar@<-2pt>[u]\ar[ll]&7\ar[l]\ar[rd] \\
&& 2\ar[ru]\ar[ld] &&& 3\ar[rru]\ar[llldd] &&&4\ar[uuuu]\ar[dd]\\
1\ar@<-2pt>[r]\ar[rruuuuu] & 14\ar@<-2pt>[l]\ar[rd] &&&&&& && 12\ar[lu]\ar@<-2pt>[r] &5\ar[llllldddd]\ar@<-2pt>[l] \\
&&20\ar[uu]\ar[dddrrr] &&&&&&6\ar[llluu]\ar[ru] \\
\\\\
&&&&& 13\ar[llllluuuu]\ar[rrruuu]
}\]
\[\xymatrix@R=19pt@C=19pt{
&&     &&&&&&  \\
&&& \kb \ar[rd]^{
\left[\begin{smallmatrix}
 1\\0
\end{smallmatrix}\right]
} &&&&  \kb \ar[dd]^1\ar[llll]_1\\
&&&&  \kb ^2  \ar@<-2pt>[d]_{\left[\begin{smallmatrix}
 0&1
\end{smallmatrix}\right]} \ar[rr]^1&&  \kb ^2 \ar@<-2pt>[d]_{\left[\begin{smallmatrix}
 0&1
\end{smallmatrix}\right]} \ar[ru]^{
\left[\begin{smallmatrix}
 0&1
\end{smallmatrix}\right]
}\\
&&&  \kb \ar[uu]^1\ar[rrd]_0& \kb \ar[l]_1 \ar@<-2pt>[u]_{\left[\begin{smallmatrix}
 1\\0
\end{smallmatrix}\right]}&& \kb\ar@<-2pt>[u]_{\left[\begin{smallmatrix}
 1\\0
\end{smallmatrix}\right]} \ar[ll]_1& \kb \ar[l]_1\ar[rd]^1 \\
&&  \kb \ar[ru]^1 &&&  \kb \ar[rru]_1\ar[llldd]^1 &&& \kb \ar[dd]_0\\
&  &&&&&& &&  \kb \ar[lu]_1 & \\
&& \kb \ar[uu] ^1&&&&&& \kb \ar[llluu]^1\ar[ru]_1 \\
\\\\
&&&&
}\]
\caption{The quiver $Q$ and the representation $T(1)$ for the knot diagram of Figure~\ref{fig knot 10}.}
\label{fig 515}
\end{center}
\end{figure}

\section{Kauffman states and submodules of the link diagram module}\label{sect 3}

We keep the notation of the previous sections. Again we choose a segment $i\in K_1$ and consider the Kauffman states and the $B$-module $T(i)$. 
Our goal is now to prove that the lattice of Kauffman states of a  link $K$ relative to a segment $i$ is isomorphic to the lattice of submodules of the direct summand $T(i)$ of the corresponding  link diagram  module $T$.

\subsection{The state module $M(\S)$} 
Let $\S$ be a Kauffman state. We will define a $B$-module $M(\S) =(M_j,M_\za)_{j\in Q_0,\za\in Q_1}$. Consider a sequence $\s$ of counterclockwise transpositions that transforms the minimal Kauffman state into the state $\S$. Then we define \[M_j=\kb^{e_j},\]
where $e_j$ is the number of occurrences of $j$ in $\s$.
 The order in which the transpositions at $j$ occur determines a basis for $M_j$, which we call the \emph{basis induced by} $\s$. 

Next we define the linear maps of the representation $M$.  In the remainder of this section, we use the following matrices
\[ J_\ell =
  \left(\begin{array}{c|ccc}
    0  & \multicolumn{3}{c}{\multirow{3}{*}{{\scalebox{1.5}{$I_{\ell}$}}}} \\
    \raisebox{0pt}{\vdots} \\0 & & &\\ \hline
    0 &0 &\cdots & 0
  \end{array}\right) 
  \quad
  V_\ell =
  \left(\begin{array}{c|ccc}
       0  & \multicolumn{3}{c}{\multirow{3}{*}{{\scalebox{1.5}{$\ I_\ell\ $}}}} \\
    \raisebox{2pt}{\vdots} & & &\\
    0 & & & 
  \end{array}\right)
\quad
H_\ell=  \left(\begin{array}{cccc}
  \multicolumn{3}{c}{\multirow{3}{*}{\raisebox{0mm}{\scalebox{1.5}{$I_{\ell}$}}}} 
  \\ 
    \\
    \\ \hline
0&\cdots&0 
  \end{array}\right)
\]
where $I_\ell$ denotes the identity matrix of size $\ell$.  We point out that $J_\ell$ is a Jordan block of size $\ell+1$ with eigenvalue 0, and that $H_\ell V_\ell=J_\ell$ and $V_\ell H_\ell = J_{\ell -1}$.

Every crossing point $p\in K_0$ indues a subsequence $\s(p)$ of $\s$ consisting of all occurrences of the transpositions at the four segments incident to $p$.
Let $a,b,c,d$ denote these four segments in counterclockwise order around $p$ such that $a$ is the first entry in $\s(p)$. Then $\s(p)$ is of one of the following forms
\begin{equation}
\label{eq wp} (abcd)^\ell,\ (abcd)^\ell a,\ (abcd)^\ell ab,\ (abcd)^\ell abc,
\end{equation}
for some $\ell\ge 0$. Let 
\[ 
\xymatrix{a\ar[r]^\zd&d\ar[d]^\zg\\
b\ar[u]^\za&c\ar[l]^\zb}\] 
be the corresponding 4-cycle in the quiver $Q$, and let $w_p=\zd\zg\zb\za\in B$.

Since every arrow of $Q$ lies in a unique 4-cycle induced by a crossing point, it suffices to define the linear maps of the representation $M$ on these four arrows $\za,\zb,\zg,\zd$. There are four cases depending on the sequence $\s(p)$.
\begin{itemize}
\item [(i)] 
If $\s(p)=(abcd)^\ell$ then $e_a=e_b=e_c=e_d=\ell$ and
\[
M_\zd= J_{\ell-1}
 \quad M_\zg=M_\zb=M_\za= I_\ell.\]

\item [(ii)] If $\s(p)=(abcd)^\ell a$ then $e_a=\ell+1, e_b=e_c=e_d=\ell$ and
\[
M_\zd=V_{\ell}
 \quad M_\zg=M_\zb= I_\ell, \quad M_\za= 
H_{\ell} 
\]

\item [(iii)] If $\s(p)=(abcd)^\ell ab$ then $e_a=e_b=\ell+1, e_c=e_d=\ell$ and
\[
M_\zd=V_{\ell}
 \quad M_\zg= I_\ell, \quad M_\zb= 
 H_{\ell} \quad M_\za=I_{\ell+1} .\]

\item [(iv)] If $\s(p)=(abcd)^\ell abc$ then $e_a=e_b=e_c=\ell+1, e_d=\ell$ and
\[
M_\zd=V_{\ell}
 \quad M_\zg=H_{\ell}
 \quad M_\zb=M_\za=I_{\ell+1} .\]
\end{itemize}
 
\begin{definition}\label{def M(S)}
 The $B$-module $M(\S)$ is called the \emph{state module} associated to the Kauffman state $\S$.
\end{definition}
\begin{remark}\label{rem::wp}
 In all four cases (i)-(iv) above the composition of the four matrices along the cycle $w_p$ is equal to $J_{\ell-1}$. Thus the action of $w_p$  on $M(\S)$ is given by this matrix.
\end{remark}

From the construction of the state module, we have the following results.

\begin{lemma}
 \label{lem e}
 Let $M(\S)=(M_x,M_\za)_{x\in Q_0, \za\in Q_1}$. Then for every arrow $\za\colon j\to k$, we have $|\dim M_j-\dim M_k|\le 1$.
\end{lemma}

\begin{lemma}
 \label{lem e2} If the state $\S'$ is obtained from the state $\S$ by applying the transposition at a segment $a$ then the module $M(\S')$ is obtained from $M(\S)$ by
 \begin{itemize}
\item [(i)] increasing the dimension at vertex $a$ by one;
\item [(ii)] increasing the rank of the map on  each arrow $\za\colon a\to \bullet$ starting at $a$ by one.
\end{itemize}
The dimension at the other vertices and the rank on the other arrows do not change.
\end{lemma}

We also note the following for future reference. 
\begin{lemma}\label{lem 5.4}
 Let $\S$ be a state and $\s$ a sequence of transpositions that transforms the minimal state into $\S$. If \[w=\xymatrix{a_0\ar[r]^{\za_1}& a_1\ar[r]^{\za_2}&\dots\ar[r]^{\za_{t-1}}& a_t\ar[r]^{\za_t}& a_0}\] 
 is a chordless cycle in $Q$ then the subsequence of $\s$ of all occurrences of transpositions at vertices of $w$ is of the form
 \[ a_j\dots a_2 a_1 a_0(a_t a_{t-1}\dots a_1 a_0)^\ell a_t a_{t-1}\dots a_k\] for some $j,k$ and $\ell$. In particular, the order in $\s$ is opposite to the order in $w$.
\end{lemma}
\begin{proof} We have already proved this result in equation~(\ref{eq wp})  in the case where $w=w_p$ is the chordless 4-cycle given by a crossing point $p\in K_0$. It suffices to show the result in the case where $w=w_R$ is the chordless cycle of a region $R\in K_2$.
 The transposition at $a_i$ is defined by moving the markers counterclockwise around the  endpoints of $a_i$, but it can also be seen as moving the markers in the clockwise direction along the segment $a_i$, see 
Figure~\ref{fig::KauffmanTransposition}. Thus the proof for the region cycle $w_R$ is dual to the proof for the crossing point cycle $w_p$.
%\begin{figure}
%\begin{center}\scalebox{1}{\input{figtranspositionregion.pdf_tex}}
%\caption{The  transposition at segment $a_i$ transforms the picture on the left to the one on the right. It can be see as moving the markers counterclockwise around the crossing points into the next region, or it can be see as moving the markers clockwise around the regions up to the next crossing point.}
%\label{fig transposition region}
%\end{center}
%\end{figure}
\end{proof}
As an immediate consequence we have the following.
\begin{corollary}\label{cor e}
 If $\za\colon a\to d$ is an arrow in $Q$ such that the transpositions at $a$ and $d$ occur consecutively in $\s$ then $d$ comes before $a$. \qed
\end{corollary}

\subsection{Lattice isomorphism}
We start with the following result. 
\begin{proposition}
 \label{prop a-f}
Let $K$ be a  link diagram without curls and let $i\in K_1$ be a segment. Let $\S,\S'$ be two Kauffman states relative to $i$. Then
\begin{itemize}
\item [(a)] $M(\S)$ is a $B$-module.
\item [(b)] $M(\S)\not\cong M(\S')$ if $\S\ne \S'$.
\item [(c)] If $\S$ is the minimal Kauffman state, then $M(\S)=0$.
\item [(d)] If $\S$ is the maximal Kauffman state, then $M(\S)=T(i)$.
\item [(e)] If $\S<\S'$ then $M(\S)$ is a submodule of $M(\S')$.
\item [(f)] For every submodule $M$ of $T(i)$ there is a unique Kauffman state $\S$ such that $M\cong M(\S)$.
\end{itemize}
\end{proposition}

\begin{proof}
 (a) By definition, $M(\S)$ is a representation of $Q$, so we only need to check that $M(\S)$ satisfies the relations given by the cyclic derivatives of the potential $W$. 
 
 Let $(p,R)$ be an arrow in $Q$, thus $p\in K_0$ and $R\in K_2$ such that the region $R$ is incident to the crossing point $p$. By definition of the potential, we have 
 \[\partial_{(p,R)} W = w_p-w_R,\]
 where $w_p=(p,R) w'_p$ is the 4-cycle in $Q$ given by the four arrows around the crossing point $p$ and
 $w_R=(p,R)w_R'$ is the cycle in $Q$ given by the arrows around the region $R$, see Figure~\ref{fig point region}. \begin{figure}
\begin{center}
\begin{minipage}{2in} 
%% Creator: Inkscape 1.0 (4035a4f, 2020-05-01), www.inkscape.org
%% Accompanies image file 'figpointregion.pdf' (pdf, eps, ps)
\begingroup%
  \makeatletter%
  \providecommand\color[2][]{%
    \errmessage{(Inkscape) Color is used for the text in Inkscape, but the package 'color.sty' is not loaded}%
    \renewcommand\color[2][]{}%
  }%
  \providecommand\transparent[1]{%
    \errmessage{(Inkscape) Transparency is used (non-zero) for the text in Inkscape, but the package 'transparent.sty' is not loaded}%
    \renewcommand\transparent[1]{}%
  }%
  \providecommand\rotatebox[2]{#2}%
  \newcommand*\fsize{\dimexpr\f@size pt\relax}%
  \newcommand*\lineheight[1]{\fontsize{\fsize}{#1\fsize}\selectfont}%
  \ifx\svgwidth\undefined%
    \setlength{\unitlength}{99.76402255bp}%
    \ifx\svgscale\undefined%
      \relax%
    \else%
      \setlength{\unitlength}{\unitlength * \real{\svgscale}}%
    \fi%
  \else%
    \setlength{\unitlength}{\svgwidth}%
  \fi%
  \global\let\svgwidth\undefined%
  \global\let\svgscale\undefined%
  \makeatother%
  \begin{picture}(1,0.60850702)%
    \lineheight{1}%
    \setlength\tabcolsep{0pt}%
    \put(0,0){\includegraphics[width=\unitlength,page=1]{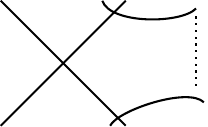}}%
    \put(0.27418253,0.34935992){\makebox(0,0)[lt]{\lineheight{1.25}\smash{\begin{tabular}[t]{l}$p$\end{tabular}}}}%
    \put(0.42453737,0.19148739){\makebox(0,0)[lt]{\lineheight{1.25}\smash{\begin{tabular}[t]{l}$a_k$\end{tabular}}}}%
    \put(0.4228581,0.36842348){\makebox(0,0)[lt]{\lineheight{1.25}\smash{\begin{tabular}[t]{l}$a_{k+1}$\end{tabular}}}}%
    \put(0.65758722,0.289218){\makebox(0,0)[lt]{\lineheight{1.25}\smash{\begin{tabular}[t]{l}$R$\end{tabular}}}}%
  \end{picture}%
\endgroup%

\end{minipage}
\begin{minipage}{2in}
\xymatrix{%&&\cdot\ar[ld]\\
\cdot\ar[r]&a_{k+1}\ar[d]^{(p,R)}&\cdot\ar[l] \\
\cdot\ar[u]&a_k\ar[l]\ar[r]&\cdot\ar@{.}[u]
}
\end{minipage}
\caption{A crossing point $p$ with an adjacent region $R$ in the  link diagram on the left and the corresponding paths in the quiver on the right. The cycle $w_p$ is the 4-cycle. }
\label{fig point region}
\end{center}
\end{figure}
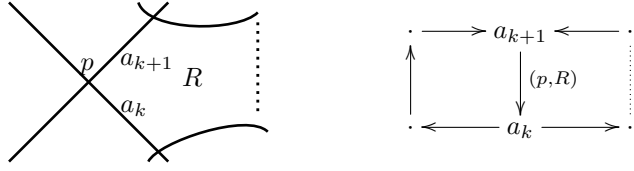
We must show that $\partial_{(p,R)} W$ acts trivially on $M(\S)$, and for that it suffices to show that the composition of the linear maps in the representation $M(\S)$ along the paths $w'_p$ and $w'_R$ are equal.

We will write $M_w$ for the composition of the arrows in $M(\S)$ along a path $w$.
According to Remark~\ref{rem::wp}, we have
$M_{w_p} = J_\ell$ for some $\ell\ge 0$, and thus it suffices to show $M_{w_R}=J_\ell$.

As before, we let $\s$ be the sequence of counterclockwise transpositions that transforms the minimal Kauffman state into the state $\S$. Let $\s(R)$ be the subsequence of $\s$ consisting of all occurrences of the transpositions at the segments that bound the region $R$. Let $a_0,a_1,\ldots, a_t$ denote these segments in clockwise order around $R$ such that $a_0$ is the first entry in the sequence $\s(R)$. By Lemma~\ref{lem 5.4}, the subsequence $\s(R)$ is of the form 
\[ \s(R)=(a_0 a_1\ldots a_t)^\ell a_0a_1\dots a_u , \ \textup{with $u<t$ or }\]
\[ \s(R)=(a_0 a_1\ldots a_t)^\ell \]
for some $\ell\ge 0$.
In the first case, the dimension of $M(\S)$ is $\ell$ at vertices $a_{u+1},\dots ,a_t$ and it is $\ell+1$ at vertices $a_0,a_1,\dots,a_u$,
while in the second case, the dimension is $\ell$ at all vertices $a_0,a_1\dots ,a_t$.

Denote the crossing point of the segments $a_{j-1}$, $a_j$ by $p_j$. So the arrow $(p_j,R)$ is $a_j\to a_{j-1}$. Then by definition of $M(\S)$ at the crossing point $p_j$, we have
\[M_{(p_j,R)} = J_{\ell-1} \quad \textup{ if $j=0$ and $\s(R)=(a_0,\dots,a_t)^\ell$}\]
and otherwise
\[M_{(p_j,R)} = \left\{
\begin{array}
 {ll}
 I_{d} & \textup{if $\dim M_{a_j}=\dim M_{a_{j-1}}=d$;}\\
 V_{d} & \textup{if $\dim M_{a_j}-1=\dim M_{a_{j-1}} =d$;}\\
 H_{d}
&\textup{if  $\dim M_{a_j}=\dim M_{a_{j-1}}-1=d$.}\\
\end{array}\right.\]

In particular
$M_{w_R}=J_\ell$.
This shows that $M(\S)$ satisfies all relations of the form $\partial_\za W$, $\za\in Q_1$.

We also have to consider the closure $I$ of the ideal generated by the relations $\partial_\za W$. For this we must show that arbitrary long paths act as zero; more precisely, if $w$ is a path such that for all $N$ there exists  a path $u_n$ of length $n>N$ such that $w=u_n$ then $M_w=0$. Suppose $w$ is such a path. Then,
for all $m$ there exists $n$ such that there is an $x\in Q_0$ through which $u_n$ passes at least $m$ times. Thus $u_n $ decomposes as 
\[ u_n= u w_1 w_2\dots w_m v,\]
where each $w_i$ is an oriented cycle that starts and ends at $x$ and that does not pass through $x$ another time. 
Take $m>\dim M_x=d$. We shall show below that, on each cycle $w_i$, the matrix product $M_{w_i}$ is some power of the matrix $J_{d-1}$. Therefore $M_{w_1w_2\dots w_n}=(J_{d-1})^{m+k}$, which is zero. Thus $M_{u_n}=0$. Since $u_n=w$, we have 
$M_w=0$, as desired. 

It remains to show that, for every oriented cycle \[w=\xymatrix{a_0\ar[r]^{\za_1}& a_1\ar[r]^{\za_2}&\dots\ar[r]^{\za_{t-1}}& a_t\ar[r]^{\za_t}& a_0}\] in $Q$ such that $a_i\ne a_j$ if $i\ne j$, the matrix  $M_w$ is a power of the matrix $J_\ell$, for some $\ell$.
Let $w$ be such a cycle. %We want to show that $M_w=J_{\ell}^m$, for some $m\ge 1$.
  By definition of $M(\S)$, for every arrow $\za_j$, the matrix $M_{\za_j}$ is one of the four matrices $I,V,H,J$, which satisfy the relations $H_\ell V_\ell=J_\ell$ and $ V_\ell H_\ell=J_{\ell-1}$. Hence if $M_w$ is not of the claimed form then $M_w=I_\ell$ and all $M_{\za_j}=I_{\ell}$. Then $\dim M_{a_j}=\ell$ at each vertex $a_j$ in $w$, which means that the transposition at $a_j$ appears exactly $\ell$ times in the sequence $\s$. Let $a_k$ be the last transposition in the sequence $\s$ at a vertex in $w$. Consider the crossing point $p$ where $a_{k-1}$ and  $a_k$ meet in the  link diagram. By Lemma~\ref{lem e2}, the last transposition at $a_k$ does not increase the rank of the matrix $M_{\za_k}$, since the arrow $\za_k$ ends in $a_k$. Thus $M_{\za_k}=J_{\ell-1}$ and hence $M_w\ne I_\ell$, and we are done.
This completes the proof of part (a) of the proposition.

%In a similar way, we can show that $M_w=J_{\ell-1}$ for some $\ell$, for every oriented cycle $w$ that does not pass twice through the same vertex, except that its starting point is equal to its endpoint. Indeed, let $\s(w)$ be the subsequence of $\s$ consisting of all occurrences of transpositions along $w$. Let $w=a_0\to a_1\to\dots\to a_t\to a_{t+1}=a_0$. Then by the same argument as above, $\s(w)=(a_0,\dots,a_t)^\ell a_0,\dots,a_u$ for some $\ell>0$ and $ u<t$, or
%$\s(w)=(a_0,\dots,a_t)^\ell.$
%By definition of $T(i)$, for every arrow $\za_j$, the matrix $M_{a_j}$ is one of the four matrices $O,V,H,I$. 
%
% Again we conclude $M_w=J_{\ell-1}$.

(b) Let $\S,\S'$ be two Kauffman states and suppose that $M(\S)=M(\S')$. Let $\s$ and $\s'$ be the sequences of transpositions that transform the minimal state into the the state $\S$ and  $\S'$, respectively.  Let $p$ be any crossing point and denote by $\s(p)$ and $\s'(p)$ the subsequences of $\s$ and $\s'$ consisting of all occurrences of transpositions at $p$. Since $M(\S)\cong M(\S')$, both representations have the same dimension vector and thus $\s(p)$ and $\s'(p)$ are equal up to a permutation. In fact,  since the minimal state has exactly one marker at the point $p$, it follows that $\s(p)=\s(p')$. 
At every crossing point $p$, the states $\S$ and $\S'$ are determined by the last entry in $\s(p)=\s'(p)$, and thus $\S=\S'$.

(c) If $\S$ is the minimal state then its sequence of transpositions $\s$ is empty. Thus $M(\S) $ is the zero module.

(d) Let $\S$ be the maximal state. 
We have described Kauffman's construction of $\S$ in \cite{K} as a partition of $K_1$ in section~\ref{sect partition}. The fact that $M(\S)$ and $T(i)$ have the same dimension vector follows directly from that. 
 We now show that $M(\S)$ and $T(i)$ also have the same linear maps. 

Let $\za\colon a\to d$ be an arrow in $Q$. It is clear from the definition of $M(\S)$ and $T(i)$ and by Lemma~\ref{lem e} that the linear maps on $\za$ are the same if the dimension at vertex $a$ is different from the dimension at vertex $d$.  Suppose therefore that $\dim M_a=\dim M_d=\ell$. Recall that the arrow $\za$ corresponds to a pair $(p,R)$, where $p$ is a crossing point and $R$ is an adjacent region in $K$. In the quiver $Q$, we have two corresponding chordless cycles $w_p$ and $w_R$ that share the arrow $\za$ as follows.
\[\xymatrix{b\ar[r]&a\ar[d]^\za&a_1\ar[l] \\
c\ar[u] &d\ar[l]\ar[r]&a_{t-1} \ar@{.>}[u]}\]
The crossing point cycle $w_p$ is the cycle of length 4 on the left. The length of the region cycle $w_R$ is the number of segments that bound $R$ in the  link diagram. We denote this length by $t+1$. 
Note that in these two cycles the arrow $\za $ is the only arrow that ends at $d$.  Consider the sequence of transpositions $\s$ that transforms the minimal state into the maximal state $\S$ and let $\s(p)$ and $\s(R)$ be the subsequences of all occurrences of transpositions at segments incident to $p$, respectively at segments bounding $R$. Suppose first that the first occurrence of $a$ is before the first occurrence of $d$ in $\s(p)$. Then $a$ must be the first entry in $\s(p)$, because the direct predecessor of $a$ would have to be $d$, by Corollary~\ref{cor e}. Similarly, $d$ must be the last entry in $\s(p)$, because $M_a$ and $M_d$ have the same dimension $\ell$ and the direct successor of $d$ in $\s(p)$ would have to be $a$.  
Thus $\s(p)=(abcd)^\ell$ and this shows that
\begin{itemize}
\item[(i)] $M_\za=J_{\ell -1}$, by definition of $M(\S)$; 
\item[(ii)] The dimension of $M(\S)$ is $\ell$ at each vertex of $w_p$ and thus $p$ is an internal point of $K'(\ell)$ as in Definition~\ref{def R(x)}. 
\end{itemize}

Now consider the sequence $\s(R)$. Since $a$ occurs before $d$ in $\s(p)$, it also does so in $\s(R)$, and by the same argument as above, we see that $a$ must be the first entry of $\s(R)$ and $d$ must be the last. Thus $\s(R)=(aa_1a_2\dots a_{t-1} d)^\ell $
and therefore the dimension of $M(\S)$ is equal to $\ell$ at every vertex of $w_R$.  
The transposition at $a$ moves  two state markers at the endpoints of $a$ counterclockwise. By our convention on the orientation of the quiver,  the marker at the endpoint $p$ must lie in the region $R$. 
Moreover, the fact that $a$ is the first entry in both sequences $\s(p)$ and $\s(R)$ implies that the position $(p,R)$ carries a state marker already in the minimal state. Similarly, since $d$ is the last entry in $\s(p)$ and $\s(R)$, the position $(p,R)$ also has a state marker in the maximal state. It follows from the construction of the minimal and maximal states in \cite{K} that the region $R$ is the region $R(p)$ of the internal point $p$ as in Definition~\ref{def R(x)}. Now the definition of the maps in
the   diagram module $T(i)$
implies that $T(i)_\za=J_{\ell -1}$. Hence $M(\S)_\za=T(i)_\za$. 

It remains the case where the first occurrence of $d$ is before the first occurrence of $a$ in $\s(p)$. Then Lemma~\ref{lem e2} implies that each occurrence of $a$ in $\s(p)$ augments the rank of $M(\S)_\za$ by one. Thus our assumption $\dim M_a=\dim M_d=\ell$ implies that $M_\za=I_\ell$. On the other hand, we also have $T(i)_\za=I_\ell$, because the position $(p,R)$ does not carry the state marker of maximal state, and thus $R$ is not the region $R(p)$ of the internal point $p$. This completes the proof of part (d). 

(e) It suffices to show that if the state $\S'$ is obtained from the state $\S$ by a single transposition at some segment $a$ then $M(\S)$ is a submodule of $M(\S')$.  
 We use the notation $M(\S)=(M_x,M_\zg)$, $M(\S')=(M_x',M_\zg')$ and $d_x=\dim M_x$, $d_x'=\dim M_x'$. Define a morphism $f\colon M(\S)\to M(\S')$ by
\[f_j=\left\{
\begin{array}{ll}
I_{d_j}&\textup{if $j\ne a$};\\
H_{d_j}&\textup{if $j=a$.}
\end{array}\right.
\]
Clearly $f$ is injective. To show that $f$ is a morphism of $B$-modules, we need to consider arrows $\xymatrix{b\ar[r]^\za&a\ar[r]^\zb&c}$ in $Q$ and show that the following diagram commutes.
\[\xymatrix{M_b\ar[r]^{M_\za}\ar[d]^{f_b}  & M_a\ar[r]^{M_\zb}\ar[d]^{f_a}   
 & M_c\ar[d]^{f_c}   
\\
M_{b}'\ar[r]^{M_\za'}&
M_{a}'\ar[r]^{M_\zb'}&
M_{c}'}\]
Since $\S'$ is obtained from $\S$ by the transposition at $a$, Lemmata~\ref{lem e} and \ref{lem e2} imply that $d_a\in\{d_b,d_b-1\}=\{d_c,d_c -1\}$. Moreover $d_a'=d_a+1$, $d_b'=d_b$ and $d_c'=d_c$, and the maps $M_\za, M_\zb,M_\za', M_\zb'$ are uniquely determined by the fact given in Lemma~\ref{lem e2} that 
$\textup{rank}(M_\za')=\textup{rank}(M_\za)$ and
$\textup{rank}(M_\zb')=\textup{rank}(M_\zb)+1.$
 Thus
 \[ 
\begin{array}
 {llll}
 M_\za =I_\ell &\textup{and}&M_\za'=H_\ell &\textup{if $d_a=d_b=\ell$;}\\
 M_\za =V_\ell &\textup{and}&M_\za'=J_\ell &\textup{if $d_a=d_b-1=\ell$;}\\
 M_\zb =J_{\ell-1} &\textup{and}&M_\zb'=V_\ell &\textup{if $d_a=d_c=\ell$;}\\
 M_\zb =H_\ell &\textup{and}&M_\zb'=I_{\ell+1} &\textup{if $d_a=d_c-1=\ell$.}\\
 \end{array}\]
 This shows that the diagram commutes and the proof of (e) is complete.

(f) Let $\S_{max}$ denote the maximal Kauffman state. Thus $T(i)=M(\S_{max})$ by part (d).  Let $d_a=\dim T(i)_a$.
We fix a sequence of transpositions $\s$ that transforms the minimal state into the maximal state, and we use the basis of $T(i)$ induced by $\s$. In particular, we have a basis $\{e_1,e_2,\dots,e_{d_a}\}$ for every vector space $T(i)_a$ with $a\in Q_0$.

Let $M=(M_a,M_\za)$ be a submodule of $T(i)$. Each vector space $M_a$ is a subspace of $T(i)_a$, and thus the points in $M_a$ can be expressed as coordinate vectors $(x_1,\dots,x_{d_a})$ with respect to our basis of $T(i)_a$. For every vertex $a\in Q_0$, let $\pi_a$  denote the canonical projection from the vector space $M=\oplus_{j\in Q_0} M_j$ to the vector space $M_a$. For any point $x\in M$, we define an integer $m(a,x)$ as follows. If $\pi_a(x)=(x_1,\dots,x_{d_a})\ne 0$, we let $m(a,x)$ be the unique integer such that $x_{m(a,x)}\ne 0$ and $x_{k}=0$, for all $k=m(a,x)+1,\dots ,d_a$. If $\pi_a(x)=0$, we let $m(a,x)=0$. We then define a function 
$m\colon Q_0\to \Z$ by $m(a)=\max_{x\in M} m(a,x)$.
We will show that $m(a)=\dim M_a$.

If $m(a)=0$ then $\pi_a(x)=0$, for all $x\in M$, and thus $M$ is not supported at $a$, whence $\dim M_a =0.$
If $m(a)= 1$ then $\pi_a(x)\in \textup{span}\{e_1\}$, for all $x\in M$, and hence $\dim M_a=1$.  

Suppose now that $m(a)\ge 2$. Then $\dim T(i)_a \ge 2$ and thus the sequence $\s$ contains the transposition at $a$ at least twice. Denote the crossing points at the ends of the segment $a$ in $K$ by $p$ and $p'$, and denote the adjacent segments by $b,c,d$ and $b',c',d'$ in counterclockwise order as shown in the left picture of Figure~\ref{fig doublecross}.
\begin{figure}
\begin{center}
\begin{minipage}
{2in} 
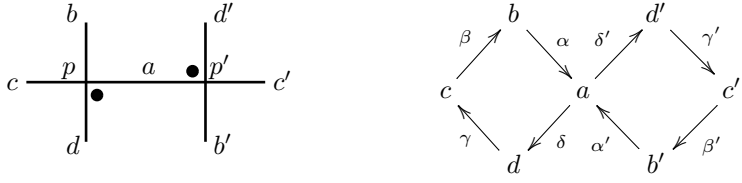
\end{minipage}
\begin{minipage}
 {2in}
 \[\xymatrix@R15pt@C15pt{&b\ar[rd]^\za&&d'\ar[rd]^{\zg'}\\
 c\ar[ru]^\zb&&a\ar[ru]^{\zd'}\ar[ld]^\zd&&c'\ar[ld]^{\zb'}\\
 &d\ar[lu]^\zg&&b'\ar[lu]^{\za'}}\]
\end{minipage}
\caption{A local configuration in the  link diagram on the left and the corresponding configuration of the quiver on the right.}
\label{fig doublecross}
\end{center}
\end{figure}
The sequence $\s$ must also contain the transpositions at $b,c,d,b',c'$ and $d'$. The corresponding subquiver of $Q$ is shown on the right in Figure~\ref{fig doublecross}. We have chordless cycles $w_p=\zd\zg\zb\za$ and $w_{p'}=\zd'\zg'\zb'\za'$.  Because of Remark~\ref{rem::wp}, the action of $w_p$ on $M_\za$ is as follows. If $\pi_a(x)=(x_1,x_2,\dots,x_{m(a)},0,\dots,0)$ then 
\begin{equation}
 \label{eq pi}
% \left\{
\begin{array}{rcl}
 \pi_a(x\cdot w_p)&=&(x_2,x_3,\dots,x_{m(a)},0,\dots,0)\\
  \pi_a(x\cdot w_p^k)&=&(x_{k+1},x_{k+2},\dots,x_{m(a)},0,\dots,0)\\
   \pi_a(x\cdot w_p^{m(a)-1})&=&(x_{m(a)},0,\dots,0)
\end{array}%\right.
 \end{equation}
 Hence, since $x_{m(a)}\ne 0$, the vectors $\pi_a(x),\pi_a(x \cdot w_p),\dots,\pi_a(x \cdot w_p^{m(a)-1})$ are linearly independent vectors in $M_a$. Thus $M_a=\textup{span}\{e_1,e_2,\ldots,e_{m(a)}\}$ and hence $\dim M_a=m(a)$ as claimed.
 
 In fact, each subspace $M_a\subset T(i)_a$, and hence the submodule $M\subset T(i)$, is completely determined by $m(a)=\dim M_a$. In particular, the submodules of $T(i)$ are determined by their dimension vector. 
 
 We will now show that $M$ corresponds to a Kauffman state using induction on the total dimension $\ell=\sum_{a\in Q_0} \dim M_a$ of $M$.
 If $\ell=0$ then $M=M(\S_{min})$ is the zero module. Suppose now that $\ell\ge 1.$ Let $a\in Q_0$ be such that $S(a)$ is a direct summand of $\textup{top} \,M$.
 Recall that $\textup{top} \,M=M/\rad M=M/(M\cdot \rad B)$ is the largest semisimple quotient of $M$. In particular, we have a short exact sequence
 \begin{equation}
\label{eq ses}
\xymatrix{0\ar[r]&L\ar[r]^f&M\ar[r]^g&S(a)\ar[r]&0,}
\end{equation}
where $L$ is the kernel of $g$.
By induction, we can assume that there exists a state $\S$ such that $L=M(\S)$. We will show that $M=M(\S')$, where $\S'$ is the state obtained from $\S$ by the transposition at the segment $a$. 

Since $M$ is determined by its dimension vector, we only need to show that the state $\S$ admits the transposition at $a$. In other words, the markers at segment $a$ must be in the positions indicated on the left of Figure~\ref{fig doublecross} and corresponding to the arrows $\zd,\zd'\in Q_1$ on the right of the same figure. It suffices to show that the subsequence $\s(p)$ of $\s$ ends in the transposition at $d$ and the subsequence $\s(p')$ ends in the transposition at $d'$. We will show this for $\s(p)$ only, since the other case is symmetric.

By Lemma~\ref{lem e}, we know that 
$|\dim L_a-\dim L_d|\le 1$ and
$|\dim M_a-\dim M_d|\le 1$.
Moreover, the short exact sequence (\ref{eq ses}) implies $\dim M_a=\dim L_a +1$, and thus $\dim L_a=\dim L_d-1$ or $\dim L_a=\dim L_d$.
In the former case, the sequence $\s(p)$ must end in $d$, and we are done.

Suppose now $\dim L_a=\dim L_d$. The morphism $g$ of $(\ref{eq ses})$ gives rise to the commutative diagram
\[\xymatrix{M_b\ar[r]^{M_\za}\ar[d] &M_a\ar[d]^{g_a}\\ 
0\ar[r]&S(a)_a}\]
with $g_a\ne 0$. Thus the commutativity implies that $g_aM_\za=0$, and hence $M_\za$ is not surjective. Because of the description of the action of $w_p$ in equation~(\ref{eq pi}), the cokernel of $M_\za\circ M_\zb\circ M_\zg\circ M_\zd$ is of dimension one, and thus the cokernel of $M_\za$ is of dimension one.
Consequently the exactness in (\ref{eq ses}) implies that the map $L_\za$ in the representation $L$ is surjective, and therefore the sequence $\s(p)$ does not end in $a$, by Lemma~\ref{lem e2}(ii). 

Next we show that $\s(p) $ does not end in $b$ or $c$. Suppose first that  $\dim M_b=\dim M_a=\ell+1$ then $\dim L_b=\dim L_a+1$ and the  sequence $\s(p)$ starts with $b$ and $L_\za=V_{\ell}$. The morphism $f$ of (\ref{eq ses}) induces a commutative diagram 

\[\xymatrix{L_b\ar[r]^{L_\za=V_\ell}\ar[d]_{f_b=I_{\ell+1}} &L_a\ar[d]^{f_a=H_\ell}\\ 
M_b\ar[r]^{M_\za}&M_a}\]
whence
$M_\za=J_\ell$.

Then $w_p=\zd\zg\zb\za$ acts like $M_\za$, and we obtain $M_\zb=M_\zg=M_\zd=I$. Thus $L_\zb=L_\zg=I$. Again Lemma~\ref{lem e2}(ii) implies that the sequence $\s(p)$ does not end in $b$ or $c$. Therefore $\s(p)$ must end in $d$, and we are done.

On the other hand,
if  $\dim M_b\ne \dim M_a$ then $\dim M_b=\dim M_a-1=\ell$ and  $\dim L_b=\dim L_a=\ell$. Then the definition of $L=M(\S)$ implies that the map $L_\za\colon L_b\to L_a$ is either $I_\ell$ or $J_{\ell-1}$. In the latter case, the sequence $\s(p)$ would end in $a$, a contradiction. Thus $L_\za=I_\ell$. In view of the sequence (\ref{eq ses}), we obtain
$M_\za=H_\ell$, and from equation (\ref{eq pi}), we have $M_\za\circ M_\zb\circ M_\zg\circ M_\zd=J_{\ell}$. In particular, $M_\zb, M_\zg$ and $M_\zd$ are surjective, and using (\ref{eq ses}) again, we see that $L_\zb$ and $L_\zg$ are surjective as well, since $\zb$ and $\zg$ are not incident to the vertex $a$. Now Lemma~\ref{lem e2}(ii) implies that $\s(p)$ does not end in $b$ or $c$. This completes the proof of part (f).
\end{proof}

\begin{corollary}
 \label{cor grassmannian}
 (a) For every dimension vector $\mathbf{e}$ the quiver Grassmannian $\textup{Gr}_{\mathbf{e}}(T(i))$ is either empty or a point. 
 In particular,  the Euler characteristic
 \[\chi(\textup{Gr}_{\mathbf{e}}(T(i)))= 0 \ or\ 1.\]

(b) The $F$-polynomial of $T(i)$ is \[F_{T(i)}= \sum_{L\subset T(i)} \mathbf{y}^{\underline{\dim}\, L},\]  where the sum is over all submodules of $T(i)$ and $\mathbf{y}^{\underline{\dim}\, L}=\prod_{i=1}^{2n}y_i^{\dim L_i}$.

\end{corollary}

\begin{proof}
In the proof of part (f) of the proposition, we have seen that every submodule of $T(i)$ is determined by its dimension vector. Thus if there is a submodule of dimension vector $\mathbf{e}$ then it is unique, and $\textup{Gr}_{\mathbf{e}}(T(i))$ is a point. Otherwise it is empty. This shows (a), and (b) follows directly.
\end{proof}

\bigskip

We are ready for the main result of this section.
\begin{theorem}
 \label{thm lattice iso}
 The map $\S\mapsto M(\S)$ is a lattice isomorphism from the lattice of Kauffman states of $K$ relative to the segment $i$ to the lattice of submodules of the module $T(i)$.
\end{theorem}
\begin{proof}
 The map is well-defined by parts (a) and (d) of Proposition~\ref{prop a-f}, injective by part (b) and surjective by part (f). Part (e) implies that it is order preserving. Moreover the maximum and minimum elements correspond by parts (c) and (d).
\end{proof}

\subsection{Indecomposability}
\begin{proposition}
The $B$-module $T(i)$ is indecomposable.
\end{proposition}
\begin{proof} First we show that the support of $T(i)$ induces a connected subquiver of $Q$. 
By definition of $T(i)$, the support consists of all segments in $K_1\setminus K(0)$, which is obtained from $K$ by removing the two regions $R_1,R_2$ that are incident to the segment $i$. If $K\setminus\{R_1,R_2\}$ is disconnected then we can draw a closed curve $\zg$ in the plane that separates one connected component from the rest. Putting back $R_1$ and $R_2$, we see that $\zg$ crosses exactly two regions in $K$. This means that $K$ is a  connected sum which contradicts our assumption that $K$ is prime. Thus the support of $T(i)$ is connected.

 Suppose now $T(i)=M\oplus N$ is the direct sum of two nonzero $B$-modules. 
 Let $\za\colon j\to k$ be an arrow in $Q$ such that $T(i)_j$ and $T(i)_k$ are nonzero.
 By definition, the linear map $T(i)_\za$ is nonzero, except possibly if there is a crossing point cycle $\xymatrix@C10pt{j\ar[r]^\za&k\ar[r]& l\ar[r]& m\ar[r]&j}$ such that $T(i)$ has dimension one at vertices $j,k,l$ and $m$, in which case one of the four maps is zero and the others are the identity. Therefore,
since the support of $T(i)$ induces a connected subquiver of $Q$,
  the supports of the two submodules $M$ and $N$ cannot be disjoint. 
 
 Then there exists $j\in Q_0$ such that $M_j$ and $N_j$ are both nonzero.  Since $j$ is a segment in $K$, there exists a crossing point $p\in K_0$ that is incident to $j$. Let $w_p$ denote the corresponding 4-cycle in $Q$.
 Because of Remark~\ref{rem::wp}, the vector space $T(i)_j=M_j\oplus N_j$ has a basis $\{e_1,e_2,\dots,e_d\}$ with respect to which the  action of $w_p$  is given by the matrix $J_{d-1}$.  
 
 Now let $m\in M_j, \, n\in N_j$ be nonzero elements and denote their expansions in the basis as
 \[m=\sum_{k=1}^n \mu_k e_k \qquad 
 n=\sum_{k=1}^n \nu_k e_k,\]
 with $\mu_k,\nu_k\in \kb$. Let $k_m$ and $k_n$ be the largest indices such that $\mu_{k_m}\ne 0$ and $\nu_{k_n}\ne 0$. 
 Then $m\cdot w_p^{k_m-1}=\mu_{k_m} e_1\in M_j$ and
  $n\cdot w_p^{k_n-1}=\nu_{k_n} e_1\in N_j$,
  because $M$ and $N$ are right $B$-modules. Dividing by the scalars shows that $e_1\in M_j\cap N_j$, a contradiction to the assumption that the sum $T(i)=M\oplus N$ is direct.
\end{proof}
%%%%%%%%%%%%%%%%%%%%%%%
%%
%% SECTION
%%
%%%%%%%%%%%%%%%%%%%%%%%
\section{The main result}\label{sect main}
In this section, we prove that the Alexander polynomial of the link  is a specialization of the $F$-polynomial of any indecomposable summand of the  link diagram module.

Let $K$ be an oriented link diagram of a prime link and assume that $K$ contains no curls. Let $n$ be the number of crossing points in $K$, and let $i$ be a segment in $K$. Let $T(i)$ be the corresponding indecomposable summand of the link diagram module and $F_{T(i)} = \sum_{L\subset T(i)} \mathbf{y}^{\underline{\dim}\, L}$ its $F$-polynomial as in Corollary~\ref{cor grassmannian}(b).
%, where the sum is over all submodules of $T(i)$ and $\mathbf{y}^{\underline{\dim}\, L}=\prod_{i=1}^{2n}y_i^{\dim L_i}$.  

For $f\in\Z[y_1,\dots, y_{2n}]$ we write $f|_t$ for the specialization of $f$ at
\begin{equation}
\label{specialization}
y_j=\left\{
\begin{array}
 {ll}
 -t &\textup{if segment $j$ runs from an undercrossing to an overcrossing;}\\  
 -t^{-1} &\textup{if segment $j$ runs from an overcrossing to an undercrossing;}\\
 -1 &\textup{if segment $j$ connects two overcrossings or two undercrossings.} \end{array}
 \right.
\end{equation}
\begin{theorem}
 \label{main thm}
 The Alexander polynomial of $K$ is equal to the specialization (\ref{specialization}) of the $F$-polynomial of every indecomposable summand $T(i)$ of the link diagram module $T$.
 That is \[\zD=F_{T(i)}|_t.\]
\end{theorem}
\begin{proof}
 Let $\zD$ denote the Alexander polynomial of $K$. 
Kauffman's theorem says that 
\begin{equation}
\label{eq 6.0}
\zD\doteq \sum_{\S} \zs(\S)\,w(\S),
\end{equation}
where the sum is over all states, $\zs(\S)=\pm 1$, and $w(\S)$ is a power of $t$. The symbol $\doteq$ means that the expressions on either side are equal up to sign and up to a power of $t$. 
We denote by $\S_{min}$ the minimal state. 
Normalizing the above identity, we find
\begin{equation}
\label{eq 6.1}
\zD\doteq\sum_{\S} \frac{\zs(\S)}{\zs(\S_{min})} 
\frac{w(\S)}{w(\S_{min})}.
\end{equation}

Let $\s=j_1,j_2,\dots,j_t$ be a sequence of transpositions that transforms $\S_{min}$ into $\S$, and let $M(\S)$ denote the state module introduced in Definition~\ref{def M(S)}. 
 Then \begin{equation}
\label{eq 6.2}
\underline{\dim}\, M(\S) = \sum_{k=1}^t \mathbf{e}_{j_k},
\end{equation}
where $\mathbf{e}_{j_k}\in\Z^{2n}$ is the vector that is 1 at position $j_k$ and 0 elsewhere.

Recall that, if a state $\S'$ is obtained from a state $\S$ by a s single transposition at a segment $j\in K_1$ then $w(j)=w(\S')/w(\S)$ is independent of the particular states $\S,\S'$ and only depends on the segment $j$.
Moreover, in this situation, Lemma 2.7 of \cite{K} implies that $\zs(\S')=-\zs(\S)$. Thus 
\begin{equation}
\label{eq 6.3}
\frac{\zs(\S)}{\zs(\S_{min})} 
\frac{w(\S)}{w(\S_{min})}
= (-w(j_1)) (-w(j_2)) \dots  (-w(j_t)).
\end{equation}
From our table in Figure~\ref{fig::weightsegment} we know that 
$-w(j)$ is equal to the specialization of $y_j$ at (\ref{eq 6.0}). Thus 
\[-w(j)=y_j|_t=\mathbf{y}^{\mathbf{e}_j}|_t.\]
Therefore the right hand side of equation (\ref{eq 6.3}) is equal to $\mathbf{y}^{\underline{\dim} \,M(\S)}|_t$. 
Applying this result to the formula in (\ref{eq 6.1}), we obtain
\[\zD\doteq \sum_\S \mathbf{y}^{\underline{\dim}\, M(\S)}|_t .\]
Now Theorem~\ref{thm lattice iso} implies 
\[\zD \doteq \sum_{L\subset T(i)} \mathbf{y}^{\underline{\dim}\, L}|_t = F_{T(i)}|_t, \] where the sum is over all submodules of $T(i)$, and thus it is the specialized $F$-polynomial.
\end{proof}

\section{A special case: 2-bridge links}\label{sect 2bridge}
\begin{figure}
\begin{center}
\scalebox{0.6}  {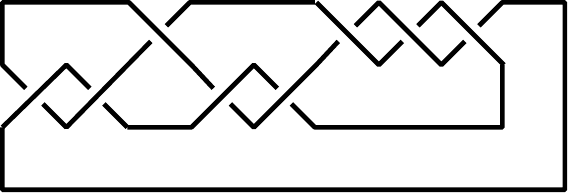}
 \caption{The 2-bridge knot $K[2,1,2,3]$. }
 \label{fig 2bridge}
\end{center}
\end{figure}

A special class of links is the family of 2-bridge links $K_{[a_1,a_2,\dots,a_n]}$ which were first studied by Schubert in \cite{Sc}. These links are parametrized by continued fractions
\begin{equation}
\label{eq cf}
[a_1,a_2,\ldots,a_n]= a_1+\cfrac{1}{a_2+\cfrac{1}{\ddots+\cfrac{1}{{a_n}}}} 
\end{equation}
with $a_i\in \mathbb{Z}_{\ge 1}$. The link $K_{[a_1,\ldots,a_n]}$ consists of $n$ braids on two strands whose number of crossings is given by the $a_i$ and that are joined together in a linear fashion as shown in Figure~\ref{fig 2bridge} and such that the resulting link is alternating. 
 $K_{[a_1,\ldots,a_n]}$  is a knot if the numerator of the continued fraction (\ref{eq cf}) is odd, and it is a link with exactly two components otherwise. For example, $K[2,1,2,3]$ in Figure~\ref{fig 2bridge} is a knot, because $ [2,1,2,3] =27/10$. 

Let $i$ be the long segment at the bottom of the link diagram that connects the $a_1$-braid to the $a_n$-braid. Then our construction of the link module $T(i)$ produces a Dynkin type $\mathbb{A}$ module whose support is given by one of the  quivers
\[\xymatrix@C15pt{1&2\ar[l] &\dots\ar[l]&\ell_1\ar[l]\ar[r] &\dots\ar[r] &\ell_2 &\dots\ar[l] &\ell_{n-1}\ar[l]\ar[r]&\dots\ar[r]&(\ell_n-1)},\]
\[\xymatrix@C15pt{1&2\ar[l] &\dots\ar[l]&\ell_1\ar[l]\ar[r] &\dots\ar[r] &\ell_2 &\dots\ar@{<-}[l] &\ell_{n-1}\ar@{<-}[l]\ar@{<-}[r]&\dots\ar@{<-}[r]&(\ell_n-1)},\]
where $\ell_j=a_1+a_2+\dots+a_j$ and direction of the arrows at the right end depend on the parity of $n$. The module $T(i)$ is of dimension one at each vertex $1,2,\dots (\ell_n-1)$ and  its linear maps on the arrows shown above are the identity maps.

In this situation, the $F$-polynomial of $T(i)$ can be computed as a sum over all perfect matchings of the snake graph $\calg_{[a_1,\dots,a_n]}$ associated to the continued fraction in \cite{CS4}. We have
\[F_{T(i)} = \sum_{P\in \textup{Match}\, \calg_{[a_1,\dots,a_n]}} y(P)
\]
where $y(P)$ is the height function of the poset of perfect matchings \cite{MS}. 

\begin{remark}
 \label{rem jones}
 This module $T(i)$ was implicitly used in \cite{LS6}, where the Jones polynomial was realized as the specialization of the $F$-polynomial of $T(i)$ at $y_1=t^{-2}$ and $y_j=-t^{-1}$ for all $j=2,3,\ldots,n$. We do not know how to generalize this specialization to other segments $i$ of this link, or to other type of links.
\end{remark}

We obtain two consequences from the above discussion.
\subsection{Type $\mathbb{A}$ cluster variables correspond to links}
We have the following result.
\begin{theorem}
 \label{thm typeA}
 Let $Q$ be a quiver of Dynkin type $\mathbb{A}$ and let $\mathcal{A}(Q)$ be its cluster algebra.  For every non-initial cluster variable $x\in \mathcal{A}(Q)$ there exists a link diagram $K$ and a segment $i\in K_1$ such that the indecomposable summand $T(i)$ of the link module is mapped to $x$ under the Caldero-Chapoton map.
\end{theorem}
\begin{proof}
 Let CC denote the Caldero-Chapoton map. Since $x$ is non-initial, there exists an indecomposable $\kb Q$-module $M$ such that $\textup{CC}(M)=x$, \cite{CCS,CC}. The support of $M$ defines a connected subquiver of $Q$ which in turn determines a continued fraction \cite{CS4} and hence a 2-bridge link $K$. From the discussion above and our main theorem, we have a segment $i\in K_1$ such that $T(i)=M$. 
\end{proof}

\subsection{An application to $q$-deformed rationals}
In \cite{MO}, Morier-Genoud and Ovsienko introduced $q$-deformed rationals and $q$-deformed continued fractions. They propose a unimodality conjecture that can be rephrased in terms of the specialized height function as follows.

Let $M$ be an indecomposable type $\mathbb{A}$ module and let $h$ be the linearization of the lattice of submodules of $M$ that  maps submodules $L$ of $M$ to their total dimension. Thus $h(L)= \dim L=\sum_{j\in Q_0} \dim L_j$.
Equivalently, we can think of $h$ as a linearization of the lattice of perfect matchings of the associated snake graph $\calg$ that maps a perfect matching $P$ to the specialization of the height function setting all $y$-variables equal to $t$. Thus
$h(P)=y(P)|_{y_j=t}$. In other words, $h$ associates to each lattice element the length of the shortest chain from the element to the minimal element in the lattice.

\begin{conjecture}
 [Morier-Genoud--Ovsienko] 
 The function $h$ is unimodal.\end{conjecture}
Progress towards this conjecture has been made in \cite{MSS}.

Using our main theorem and properties of the Alexander polynomial, we have the following result, which says that the alternating sum of the number of objects on each level of the poset is $-1,0,$ or 1.
%
%\begin{theorem}
% \label{thm height}
% Let $\calg$ be a snake graph and $\textup{Match} \,\calg $ be the poset of its perfect matchings. Let $h(P)$ be the linearization and denote by $h(p)|_{t=-1}$ its specialization  at $t=-1$. Then
% \[\sum_{P\in \textup{Match}\, \calg} h(p)|_{t=-1} = \left\{
%\begin{array}
% {ll} \pm 1&\textup{if $|\textup{Match} \,\calg| $ is odd;}\\
% 0&\textup{otherwise.}
%\end{array}\right.\]
%\end{theorem}

\begin{theorem}
 \label{thm height}
  Let $M$ be a module of Dynkin type $\mathbb{A}_n$ and $\mathcal{L}$ the submodule lattice of $M$. Then
  \[\sum_{L\in \mathcal{L}} (-1)^{h(L)}=\left\{
\begin{array}
 {ll} \pm 1&\textup{if $|\mathcal{L}| $ is odd;}\\
 0&\textup{if $|\mathcal{L}| $ is even.}
\end{array}\right.\]
\end{theorem}
\begin{proof}
 Combining Theorems~\ref{thm typeA} and \ref{main thm}, we see that there exists a 2-bridge link $K$ whose Alexander polynomial is the specialized $F$-polynomial of $M$. More precisely,
 \[\zD_K \doteq F_M|_t = \sum_{L\in\mathcal{L}} \mathbf{y}^{\underline{\dim} \,L}|_t.\]
 From the definition of the specialization (\ref{specialization}), we see that evaluating the above equation at $t=1$ gives
 \begin{equation}
\label{eq height}
\zD_K(1) = \sum_{L\in \mathcal{L}} \mathbf{y}^{\underline{\dim} \,L}|_{y_i=-1}.
\end{equation}
Furthermore
\[\mathbf{y}^{\underline{\dim} \,L}|_{y_i=-1} 
= \textstyle \prod_i y_i^{\dim\,L_i}|_{y_i=-1} 
= (-1)^{\sum_i \dim\,L_i} 
=(-1)^{h(L)}.
\]
Thus equation (\ref{eq height}) becomes
\[\zD_K(1)= \sum_{L\in \mathcal{L}} (-1)^{h(L)}.\]
Now the result follows from property (ii) of subsection~\ref{sect Alexander}
and the fact that $K$ is a knot if and only if the number of submodules of $M$ is odd.
\end{proof}

\section{Examples}

\begin{example} 
Consider the figure-eight knot. We use the same labelling of segments as in Example \ref{ex::kstateslattice}. The lattice of the sumbodules of the module $T(1)$ is shown at Figure \ref{fig::LatticeFigureEight1}. 

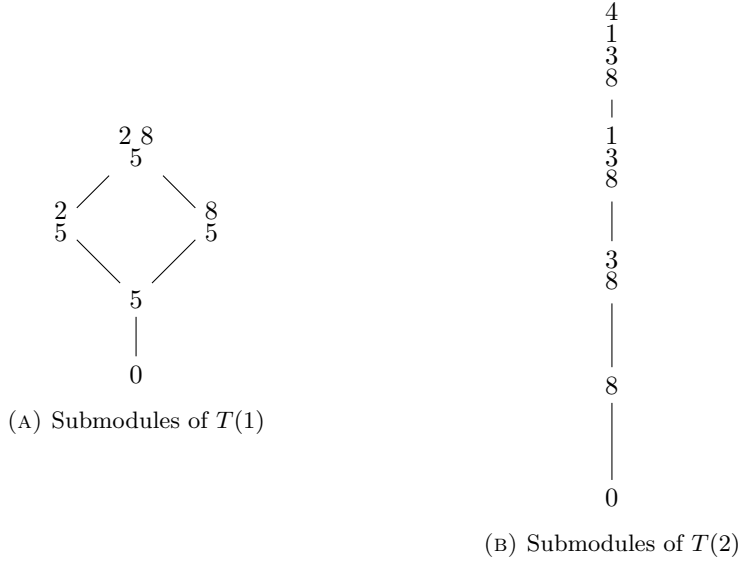
\begin{figure}[ht]
\centering
\begin{subfigure}{.5\textwidth}
    \centering
    \begin{tikzpicture}[scale = 1]
    \node (n0) at (0,0) {$\biginddeux{2 \ 8}{5}$
    };
    \node (n1) at (-1,-1) {$\biginddeux{8}{5}$
    };
    \node (n2) at (1,-1) {$\biginddeux{2}{5}$
    };
    \node (n3) at (0,-2) { $5$
    };
    \node (n4) at (0,-3) {$0$
    };
    \draw (n0) -- (n1) -- (n3) -- (n4);
    \draw (n0) -- (n2) -- (n3);
    \end{tikzpicture}
    \caption{Submodules of $T(1)$}
    \label{fig::LatticeFigureEight1}
\end{subfigure}%
\begin{subfigure}{.5\textwidth}
    \centering
    \begin{tikzpicture}[scale = 1.5]
    \node (n0) at (0,0) {
        0
    };
    \node (n1) at (0,1) {$8$
    };
    \node (n2) at (0,2) {$\biginddeux{3}{8}$
    };
    \node (n3) at (0,3) { $ \bigindtrois{1}{3}{8}$
    };
    \node (n4) at (0,4) {$\bigindquatre{4}{1}{3}{8}$
    };
    \draw (n0) -- (n1) -- (n2) -- (n3) -- (n4);
    \end{tikzpicture}
    \caption{Submodules of $T(2)$}
    \label{fig::LatticeFigureEight2}
\end{subfigure}%

\caption{Lattices of submodules of indecomposable modules on the Jacobian algebra given by the figure-eight knot.}
\label{fig::LatticeSubmoduleT(i)}
\end{figure}
 
The lattice isomporphism with the lattice of Kaufman states with regards to the segment $1$ is obvious, see Figure \ref{fig::kstateslattice}.
The $F$-polynomial of $T(1)$ is \[F_{T(1)} = 1 + y_2 + y_8 + y_2y_8 + y_2y_5y_8\] and its specialization at $y_2 = -t$, $y_5 = -t^{-1}$ and $y_8 = -t$, as given by Equation \ref{specialization}, is \[F_{T(1)}|_t = 1 -3t + t^2.\]

Remark that the lattice of sumbodules of $T(2)$ is completely different for $T(1)$, see Figures \ref{fig::LatticeFigureEight1} and \ref{fig::LatticeFigureEight2}. The $F$-polynomial of $T(2)$ is \[F_{T(2)} = 1 + y_8 + y_3y_8 + y_1y_3y_8 + y_1y_3y_4y_8\] and its specialization at $y_1 = -t^{-1}$, $y_3 = -t^{-1}$, $y_4 = -t$ and $y_8 = -t$ is \[F_{T(2)}|_t = -t^{-1} +3 -t.\]

\end{example}
\begin{example}

For the representation  $T(1)$ in Figure~\ref{fig 515}, the $F$-polynomial is the following polynomial with 75 terms. This was computed using \cite{applet}.

{\tiny\[ \begin{array}{l}
1 +y_{4 }+y_{12} y_{4 }+y_{8 }+y_{18} y_{8 }+y_{10} y_{18} y_{8 }+y_{10} y_{18} y_{19} y_{8} + 
 y_{10} y_{18} y_{19} y_{2} y_{8 }+y_{10} y_{18} y_{19} y_{2} y_{20} y_{8 }+y_{4} y_{8 }\\+y_{12} y_{4} y_{8} + 
 y_{18} y_{4} y_{8 }+y_{10} y_{18} y_{4} y_{8 }+y_{12} y_{18} y_{4} y_{8 }+y_{10} y_{12} y_{18} y_{4} y_{8} + 
 y_{10} y_{18} y_{19} y_{4} y_{8 }+y_{10} y_{12} y_{18} y_{19} y_{4} y_{8 }
 \\
 +y_{10} y_{18} y_{19} y_{2} y_{4} y_{8} + 
 y_{10} y_{12} y_{18} y_{19} y_{2} y_{4} y_{8 }+y_{10} y_{18} y_{19} y_{2} y_{20} y_{4} y_{8} 
 + 
y_{10} y_{12} y_{18} y_{19} y_{2} y_{20} y_{4} y_{8 }+y_{10} y_{18} y_{19} y_{8} y_{9} \\
+ 
 y_{10} y_{17} y_{18} y_{19} y_{8} y_{9 }+y_{10} y_{18} y_{19} y_{2} y_{8} y_{9} + 
 y_{10} y_{17} y_{18} y_{19} y_{2} y_{8} y_{9 }+y_{10} y_{18} y_{19} y_{2} y_{20} y_{8} y_{9}\\ + 
 y_{10} y_{17} y_{18} y_{19} y_{2} y_{20} y_{8} y_{9 }+y_{10} y_{18} y_{19} y_{4} y_{8} y_{9} + 
 y_{10} y_{12} y_{18} y_{19} y_{4} y_{8} y_{9 }+y_{10} y_{17} y_{18} y_{19} y_{4} y_{8} y_{9} \\+ 
 y_{10} y_{12} y_{17} y_{18} y_{19} y_{4} y_{8} y_{9 }+y_{10} y_{18} y_{19} y_{2} y_{4} y_{8} y_{9 }+ 
 y_{10} y_{12} y_{18} y_{19} y_{2} y_{4} y_{8} y_{9 }+y_{10} y_{17} y_{18} y_{19} y_{2} y_{4} y_{8} y_{9 }\\+ 
 y_{10} y_{12} y_{17} y_{18} y_{19} y_{2} y_{4} y_{8} y_{9 }+y_{10} y_{18} y_{19} y_{2} y_{20} y_{4} y_{8} y_{9 }+ 
 y_{10} y_{12} y_{18} y_{19} y_{2} y_{20} y_{4} y_{8} y_{9 }\\
 +y_{10} y_{17} y_{18} y_{19} y_{2} y_{20} y_{4} y_{8} y_{9 }+ 
 y_{10} y_{12} y_{17} y_{18} y_{19} y_{2} y_{20} y_{4} y_{8} y_{9 }+y_{10} y_{17} y_{18} y_{19} y_{4} y_{7} y_{8} y_{9 }\\+ 
 y_{10} y_{12} y_{17} y_{18} y_{19} y_{4} y_{7} y_{8} y_{9 }+y_{10} y_{16} y_{17} y_{18} y_{19} y_{4} y_{7} y_{8} y_{9 }+ 
 y_{10} y_{12} y_{16} y_{17} y_{18} y_{19} y_{4} y_{7} y_{8} y_{9 }\\
 + 
 y_{10} y_{17} y_{18} y_{19} y_{2} y_{4} y_{7} y_{8} y_{9 }+y_{10} y_{12} y_{17} y_{18} y_{19} y_{2} y_{4} y_{7} y_{8} y_{9 }+
 y_{10} y_{16} y_{17} y_{18} y_{19} y_{2} y_{4} y_{7} y_{8} y_{9 }\\+ 
 y_{10} y_{12} y_{16} y_{17} y_{18} y_{19} y_{2} y_{4} y_{7} y_{8} y_{9 }+ 
 y_{10} y_{17} y_{18} y_{19} y_{2} y_{20} y_{4} y_{7} y_{8} y_{9 }+ 
 y_{10} y_{12} y_{17} y_{18} y_{19} y_{2} y_{20} y_{4} y_{7} y_{8} y_{9 }\\+ 
 y_{10} y_{16} y_{17} y_{18} y_{19} y_{2} y_{20} y_{4} y_{7} y_{8} y_{9 }+ 
 y_{10} y_{12} y_{16} y_{17} y_{18} y_{19} y_{2} y_{20} y_{4} y_{7} y_{8} y_{9 }+ 
 y_{10} y_{17} y_{18} y_{19} y_{2} y_{20} y_{3} y_{4} y_{7} y_{8} y_{9 }\\+ 
 y_{10} y_{12} y_{17} y_{18} y_{19} y_{2} y_{20} y_{3} y_{4} y_{7} y_{8} y_{9 }+ 
 y_{10} y_{16} y_{17} y_{18} y_{19} y_{2} y_{20} y_{3} y_{4} y_{7} y_{8} y_{9 }+ 
 y_{10} y_{12} y_{16} y_{17} y_{18} y_{19} y_{2} y_{20} y_{3} y_{4} y_{7} y_{8} y_{9 }\\+ 
 y_{10} y_{12} y_{17} y_{18} y_{19} y_{2} y_{20} y_{3} y_{4} y_{6} y_{7} y_{8} y_{9 }+ 
 y_{10} y_{12} y_{16} y_{17} y_{18} y_{19} y_{2} y_{20} y_{3} y_{4} y_{6} y_{7} y_{8} y_{9 }\\+ 
 y_{10} y_{16} y_{17} y_{18} y_{19} y_{4} y_{7} y_{8}^2 y_{9 }+ 
 y_{10} y_{12} y_{16} y_{17} y_{18} y_{19} y_{4} y_{7} y_{8}^2 y_{9 }+ 
 y_{10} y_{16} y_{17} y_{18}^2 y_{19} y_{4} y_{7} y_{8}^2 y_{9 }\\
 + 
 y_{10} y_{12} y_{16} y_{17} y_{18}^2 y_{19} y_{4} y_{7} y_{8}^2 y_{9 }+ 
 y_{10} y_{16} y_{17} y_{18} y_{19} y_{2} y_{4} y_{7} y_{8}^2 y_{9 }+ 
 y_{10} y_{12} y_{16} y_{17} y_{18} y_{19} y_{2} y_{4} y_{7} y_{8}^2 y_{9 }\\+ 
 y_{10} y_{16} y_{17} y_{18}^2 y_{19} y_{2} y_{4} y_{7} y_{8}^2 y_{9 }+ 
 y_{10} y_{12} y_{16} y_{17} y_{18}^2 y_{19} y_{2} y_{4} y_{7} y_{8}^2 y_{9 }+ 
 y_{10} y_{16} y_{17} y_{18} y_{19} y_{2} y_{20} y_{4} y_{7} y_{8}^2 y_{9 }\\+ 
 y_{10} y_{12} y_{16} y_{17} y_{18} y_{19} y_{2} y_{20} y_{4} y_{7} y_{8}^2 y_{9 }+ 
 y_{10} y_{16} y_{17} y_{18}^2 y_{19} y_{2} y_{20} y_{4} y_{7} y_{8}^2 y_{9 }+ 
 y_{10} y_{12} y_{16} y_{17} y_{18}^2 y_{19} y_{2} y_{20} y_{4} y_{7} y_{8}^2 y_{9 }\\+ 
 y_{10} y_{16} y_{17} y_{18} y_{19} y_{2} y_{20} y_{3} y_{4} y_{7} y_{8}^2 y_{9 }+ 
 y_{10} y_{12} y_{16} y_{17} y_{18} y_{19} y_{2} y_{20} y_{3} y_{4} y_{7} y_{8}^2 y_{9 }+ 
 y_{10} y_{16} y_{17} y_{18}^2 y_{19} y_{2} y_{20} y_{3} y_{4} y_{7} y_{8}^2 y_{9 }\\+ 
 y_{10} y_{12} y_{16} y_{17} y_{18}^2 y_{19} y_{2} y_{20} y_{3} y_{4} y_{7} y_{8}^2 y_{9 }+ 
 y_{10} y_{12} y_{16} y_{17} y_{18} y_{19} y_{2} y_{20} y_{3} y_{4} y_{6} y_{7} y_{8}^2 y_{9 }\\+ 
 y_{10} y_{12} y_{16} y_{17} y_{18}^2 y_{19} y_{2} y_{20} y_{3} y_{4} y_{6} y_{7} y_{8}^2 y_{9}
 \end{array}
\]}
The specialization is \[F_{T(1)}|_t=3-9t+16t^2-19t^3+16t^4-9t^5+3t^6,\]
which is equal to  the Alexander polynomial of the corresponding knot $10_{66} $.
\end{example}
\begin{example}
\label{ex conway} 

The Conway knot and its quiver are illustrated at Figure \ref{fig::ConwayKnot}
\begin{figure}[ht]
\centering
        \begin{tikzpicture}[scale= 1.5]
        \path (2,1.5) (2,-.5);
        \begin{knot}[
        consider self intersections,
        %draft mode=crossings,
        clip width=5,
        ignore endpoint intersections=false,
        ]
        \strand[
            thick, 
            %show curve controls
        ] 
        (1.5,4.5) .. controls +(0,-2) and +(-0,-2) .. node[below] {$9$}
        (0,4.5) .. controls +(0,0.5) and +(0,-0.5).. node[pos=0, right] {$10$}
        (0,5) .. controls +(0,2) and +(0,2).. node[pos=0.04, right]{$11$} node[above]{$12$}
        (5,5) .. controls +(0,-0.5) and +(0.5,0).. node[pos=1, below]{$13$}
        (4,4) .. controls +(-0.5,0) and +(0.5,0).. node[pos = 1, above]{$14$}
        (3,5) .. controls +(-1,0) and +(-1,0).. node[pos = 0.4, left]{$15$} node[pos = 0.75, left]{$16$}
        (3.5,2.5) .. controls +(1,0) and +(1,0).. node[pos = 0.4, right]{$17$}
        (4,5) .. controls +(-0.5,0) and +(0.5,0).. node[pos = 0, above]{$18$}
        (3,4) .. controls +(-0.5,0) and +(0.5,0).. node[pos = 0, above]{$19$}
        (2,5) .. controls +(0,0) and +(0,0).. node[pos = 0, above]{$20$}
        (0.75,5) .. controls +(-2,0) and +(-1,1).. node[pos = 0, below]{$21$} node[pos = 0.25, above]{$22$}
        (0,2.5) .. controls +(1,-1) and +(-1,-1).. node[pos = 0.25, above]{$1$}
        (3.5,2.5) .. controls +(0.5,0.5) and +(3,0).. node[pos = 0.5, left, below]{$2$} node[pos = 0.85, below]{$3$}
        (0.75, 4) .. controls +(-2,0) and +(-1,0).. node[pos = 0, below]{$4$} node[pos = 0.25, below]{$5$} node[pos = 0.7, left]{$6$}
        (0,5.75) .. controls +(1,0) and +(0,1).. node[pos = 0.45, above, right]{$7$}
        (1.5,4.5) node[pos = 1, right]{$8$};
        \flipcrossings{1,15,4,14,7,9}
        \end{knot}
        \end{tikzpicture}
\vspace{-4cm}
        {\huge
        \[
        \scalebox{0.5}
        {\begin{xy} 0;<1pt,0pt>:<0pt,-1pt>:: 
        (130,40) *+{6} ="0",
        (250,80) *+{7} ="1",
        (250,190) *+{8} ="2",
        (210,260) *+{9} ="3",
        (170,190) *+{10} ="4",
        (170,80) *+{11} ="5",
        (450,40) *+{12} ="6",
        (410,190) *+{13} ="7",
        (370,140) *+{14} ="8",
        (330,220) *+{15} ="9",
        (370,300) *+{16} ="10",
        (450,300) *+{17} ="11",
        (410,140) *+{18} ="12",
        (370,190) *+{19} ="13",
        (290,140) *+{20} ="14",
        (210,140) *+{21} ="15",
        (130,140) *+{22} ="16",
        (130,300) *+{1} ="17",
        (370,260) *+{2} ="18",
        (290,260) *+{3} ="19",
        (210,220) *+{4} ="20",
        (130,220) *+{5} ="21",
        "5", {\ar"0"},
        "0", {\ar"6"},
        "0", {\ar"16"},
        "17", {\ar@/^50pt/"0"},
        "1", {\ar"5"},
        "6", {\ar"1"},
        "1", {\ar"14"},
        "15", {\ar"1"},
        "14", {\ar"2"},
        "2", {\ar"15"},
        "2", {\ar"19"},
        "20", {\ar"2"},
        "19", {\ar"3"},
        "3", {\ar"21"},
        "15", {\ar"4"},
        "4", {\ar"16"},
        "4", {\ar"20"},
        "21", {\ar"4"},
        "5", {\ar"15"},
        "16", {\ar"5"},
        "6", {\ar"11"},
        "12", {\ar"6"},
        "11", {\ar"7"},
        "7", {\ar"13"},
        "8", {\ar"12"},
        "14", {\ar"8"},
        "13", {\ar"9"},
        "9", {\ar"14"},
        "9", {\ar"18"},
        "19", {\ar"9"},
        "17", {\ar"10"},
        "10", {\ar"19"},
        "11", {\ar@/^30pt/"17"},
        "18", {\ar"11"},
        "16", {\ar"21"},
        "21", {\ar"17"},
        "8",{\ar@<-2pt>"13"},
        "13",{\ar@<-2pt>"8"},
        "7",{\ar@<-2pt>"12"},
        "12",{\ar@<-2pt>"7"},
        "3",{\ar@<-2pt>"20"},
        "20",{\ar@<-2pt>"3"},
        "10",{\ar@<-2pt>"18"},
        "18",{\ar@<-2pt>"10"},
        \end{xy}}\]
        } 
   
            \caption{Conway knot and its quiver}
            \label{fig::ConwayKnot}
    \end{figure}

 The $F$-polynomial of $T(1)$ has 131 terms. The highest degree term is
\[y_{2}y_{4}y_{7}y_{8}y_{10}y_{11}y_{13}y_{14}y_{15}y_{18}y_{19}y_{20}y_{21}y_{22}\]
The specialization gives $F_{T(1)}|_t=t$, confirming that the Alexander polynomial is trivial, since it is defined up to a power of $t$.
\end{example}

\pagebreak
\pagebreak
\appendix
\section{Addendum to  Knot Theory and Cluster Algebras, Adv. Math. 408  B, (2022), 108609.} 
% \author{V\'eronique Bazier-Matte}
% %\thanks{The first author was supported by }
% \address{D\'epartement de math\'ematiques et de statistique,  Universit\'e Laval, Qu\'ebec (Qu\'ebec), G1V 0A6, Canada}
% \thanks{The first author was supported by the Discovery Grants  program from The Natural Sciences and Engineering Research Council of Canada and the Research Support for New Academics from the Fonds de Recherche du Qu\'ebec Nature et Technologie}
% \email{veronique.bazier-matte.1@ulaval.ca}
% \author{Ralf Schiffler}
% \thanks{The second author was supported by the National Science Foundation grant  DMS-2054561}
% \address{Department of Mathematics, University of Connecticut, Storrs, CT 06269-1009, USA}
% \email{schiffler@math.uconn.edu}

% \maketitle
This short note is a supplement to the above article. %to our article \cite{BMS}.  
The main purpose is to improve the description of the partition $K_1=\sqcup_{d\ge 0} K(d)$ of the set of segments of a link diagram given in Section 5. % of \cite{BMS}. 
This is done in section~\ref{Addsect 1} of this note. In section~\ref{Addsect 2}, we correct errata. We thank Alfredo Najera for pointing out the first two errata to us.

\subsection{On the partition of the set of segments of a link diagram}\label{Addsect 1}

\subsubsection{Strings and their boundaries} 
\begin{figure}[htb]
\begin{center}
 \scalebox{1.5}{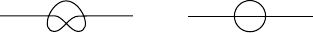}
 \caption{A trefoil string (left) and a Hopf link string (right).}
 \label{fig string}
\end{center}
\end{figure}

Let $K$ be a link diagram with $K_0, K_1,K_2 $ its sets of crossing points, segments and regions, respectively. Let $i\in K_1 $ be a segment. The \emph{$i$-th string} of $K$ is obtained from $K$ by removing an interior point from the segment $i$. The two ends of the segment $i$ that remain are called the start and the terminal of the string.  Figure~\ref{fig string} illustrates a string of the  trefoil knot and  one of the Hopf link.

When we draw the string in the plane, we can indefinitely extend its start segment horizontally to the left and its terminal segment horizontally to the right so that the complement of the string in the plane has two unbounded components. 
A segment of a string is \emph{connecting} if is incident to both unbounded components.

\begin{figure}
\begin{center}
 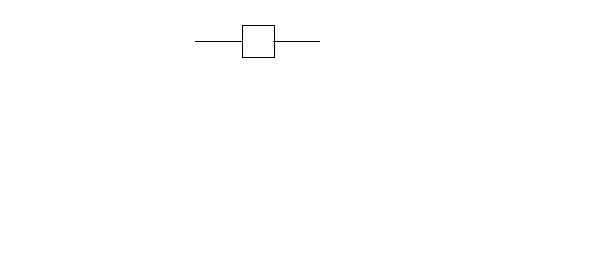
 \caption{The sum of two string $A\oplus B$ (top) and the enclosed sum of two strings $A\oplus[B,p]$ with $p$ in the interior of $A$ (center) and $p$ on the boundary of $A$ (bottom).}
 \label{fig sum}
\end{center}
\end{figure}

In \cite{K}, Kauffman defined the following operations on strings. For illustrations see Figure~\ref{fig sum}.

\begin{definition}(Sum and enclosed sum)

(a)  If $A$ and $B$ are strings, their \emph{sum} $A\oplus B$ is obtained by connecting the terminal of $A$ to the start of $B$.

(b) If $A$ and $B$ are strings, $p$ is a point on a non-connecting edge of $A$ and $p$ is not a crossing point then the \emph{enclosed sum} $A\oplus [B,p]$ of $A$ and $B$ relative to $p$ is obtained by replacing the trivial string at $p$ by $B$. In this situation, the string $A$ is called a \emph{carrier} and the string $B$ is called a \emph{rider} in $A\oplus [B,p]$.

(c) A string is called \emph{atomic} if it is irreducible with respect to both sums.

\end{definition}

\begin{remark}
    If $K$ is a prime link diagram then each string of $K$ is atomic
\end{remark}

\begin{definition} The \emph{boundary} of a string is defined as follows.
 
 (1)  If $A$ is atomic, then $\partial A$ consists of all segments of $A$ that are incident to an unbounded component. 

(2) $\partial(A\oplus B)=\partial A \oplus \partial B$.

(3) $\partial (A\oplus [B,p])=\left\{
\begin{array}
 {ll}
 \partial A \oplus [\partial B,p] &\textup{if $p\in \partial A$};\\
 \partial A &\textup{otherwise.}
\end{array}\right.
$

\end{definition}

\begin{figure}
\begin{center}
\scalebox{0.65}{ 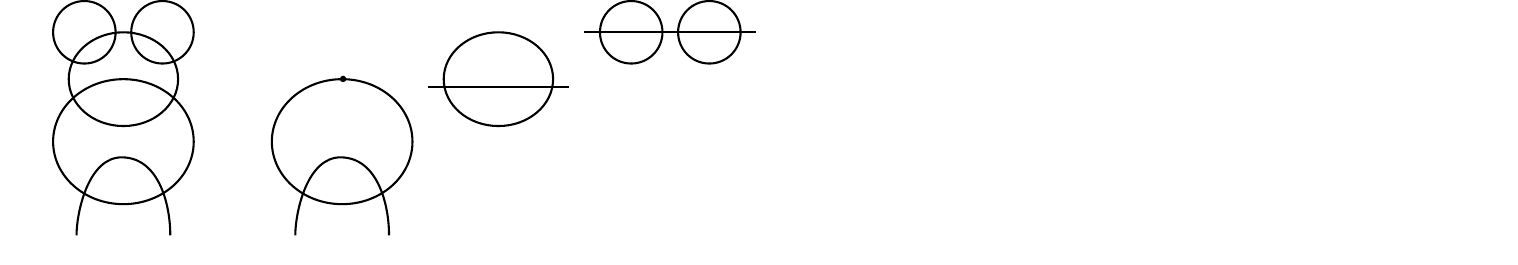}
 \caption{Computation of the boundary. The string $S$ is a double enclosed sum of the strings $A,B,C$. The boundary of $S$ is shown in the rightmost picture.}
 \label{figbdyex}
\end{center}
\end{figure}

We give an example of the computation of the boundary for a more complicated string in Figure~\ref{figbdyex}. The string there is a double enclosed sum $S=A\oplus [B \oplus[C,b],a]$ and the boundary is computed recursively  $\partial S=\partial A \oplus [\partial B\oplus[\partial C,b],a]$.

\subsection{The partition of $K_1$}
Let $i$ be a   segment of $K$. In  Section 5, we define a partition $K_1=\sqcup_{d\ge 0} K(d)$ of the set of all segments of $K$ and use it to define a representation $T(i)$. 
The sets $K(d)$ depend on the choice of the segment $i$, but, in the interest of simplicity, our notation does not reflect this dependency. This should not create confusion, since $i$ is fixed here. 

\begin{definition}
 For $d=0$, we define
\[K(0)=\{j\in K_1\mid \textup{$j$ and $i$ bound the same region of $K$}\}\cup\{i\}\]
Thus $K(0)$ is the boundary of the $i$-th string of $K$ together with $i$.

Recursively, $K(d)$ is defined as the boundary of the string (or union of strings) $K\setminus(\cup_{e<d} K(e))$.

\end{definition}

We now give an alternative description of $K(d)$, which uses the following terminology.
Given two segments $i,j\in K_1$, a curve in $\mathbb{R}^2$ is called a \emph{dimension curve} from  $j$ to  $i$ if it starts at a point on segment $j$, ends at a point on segment $i$ and does not go through a crossing point of $K$. 
Let $\dimprime(i,j)$ be the minimal number of crossings between the segments of $K$ and a dimension curve from segment $j$ to segment $i$.

With this notion, we have the following result.
\begin{proposition}
A segment $j\in K_1\setminus(\cup_{e<d} K(e))$ lies in $K(d) $ if and only if one of the following two conditions hold.
\begin{enumerate}
 \item $\dimprime(i,j)=d$, or
 \item $\dimprime(i,j)=d+1$ and $K\setminus(\cup_{e<d} K(e))$ contains a string  of the form $ A\oplus[B,p]$ for two strings $A,B$ and $p\in \partial A$ an interior point of a non-separating segment, and $j\in \partial B$ is not incident to an unbounded component of the complement of $ A\oplus[B,p]$ in the plane.
\end{enumerate}
\end{proposition}

\subsection{Errata}\label{Addsect 2}
\begin{itemize}
\item On page 10,	 in Figure 4, replace $t$ by $t^{-1}$ to match the convention elsewhere in the paper.
\item On page 16, line 3 of condition (b), rephrase as follows.
\emph{If there exists a non-constant path $w$ starting and ending at $x$ that uses only segments of $K'(d)$ and that is a subpath of the path $w_L$  or the path $ w_R$, we let $D(x)$ be the bounded domain enclosed by $w$ in the plane. }  
\item In Lemma 5.5 on page 17, the parts (b) and (c) need the additional assumption that the string $K\setminus \cup_{e<d}K(e)$ is not the sum $A \oplus B$ of two nontrivial sub-strings $A,B$.
\item Lemma  5.10 on page 19 now follows directly from the definition of boundary, because the boundary of the Hopf link string in Figure~\ref{fig string} does not contain the enclosed diameter.
\end{itemize}

% %\subsection{Examples}
% \begin{thebibliography}{}
%  \bibitem{BMS} V. Bazier-Matte and R. Schiffler,  Knot theory and cluster algebras, {\em Adv. Math.\/} {\bf 408}  B (2022), Article 108609.
% \bibitem{K} L. Kauffman, Formal Knot Theory, Mathematical Notes, 30. Princeton University Press, Princeton, NJ, 1983. 
% %
% \end{thebibliography}

\end{document}